\input amstex
\magnification=1100
\documentstyle{amsppt}
\define\al{\alpha}
\define\be{\beta}

\define\th{\theta}

\define\si{\sigma}
\define\BF{\Bbb F}
\define\BI{\Bbb I}
\define\BQ{\Bbb Q}
\define\BZ{\Bbb Z}
\define\CO{\Cal O}
\define\CR{\Cal R}
\define\gM{\goth M}

\define\bv{\bold v}
\define\bde{\bold e}
\define\oh{\overline h}
\define\ou{\overline u}

\define\ot{\overline t}
\define\oK{\overline K}
\define\oL{\overline L}
\define\oM{\overline M}

\define\tl{\tilde l}
\define\chr{\operatorname{char}}

\define\Ker{{\operatorname{Ker}}}
\define\pr{\operatorname{pr}}
\define\Cl{\operatorname{Cl}}
\define\Gal{\operatorname{Gal}}

\define\sep{{\text{sep}}}
\define\alg{{\text{alg}}}
\define\insep{{\text{insep}}}
\redefine\c{{\bold c}}
\define\cf{{\bold i}}
\redefine\i{{\bold i}}
\define\VK{V\!K}
\define\tp{^{\roman{top}}}

\define\<{\{\!\{}
\define\>{\}\!\}}
\define\isom{\simeq}

\topmatter
\author Igor Zhukov\endauthor
\title On ramification theory in the imperfect residue field case\endtitle
\address
The Department of Mathematics and Mechanics\\ St. Petersburg
University\\ Bibliotechnaya sq. 2\\ St. Petersburg 198904\\ Russia
\endaddress
\email igor@zhukov.pdmi.ras.ru \endemail
\date June 10, 2001\enddate
\thanks Supported by a grant of LMS, a research grant of Nottingham
University (ROF2/039) and grants of RFBR (projects 97-01-00058-a
and 00-01-00140)\endthanks
\endtopmatter
Let $K$ be a complete discrete valuation field with the residue
field $\oK$, $\chr \oK=p>0$. If $\oK$ is a perfect field, there
exists a beautiful theory of ramification in algebraic extensions
of $K$. Given a finite Galois extension $L/K$ with the Galois
group $G$, one can introduce a canonical filtration $(G_i)$ in $G$
with quite a natural behavior with respect to subextensions in
$L/K$. Namely, if $H$ is a normal subgroup in $G$, one has
$H_i=G_i\cap H$ and $(G/H)^j=(G^jH)/H$. In the last relation we
used upper numbering of ramification subgroups $G^j=G_{\psi(j)}$,
where the Hasse-Herbrand function $\psi=\psi_{L/K}$ can be easily
calculated in terms of orders of $G_i$. Next, this ``upper''
filtration of $G$ is compatible with class field theory. In
particular, if $L/K$ is abelian and the residue fields are finite
(or quasi-finite), then $\th(U_j)=G^j$ for all $j=0,1,\dots$,
where $\th\:K^*\to G$ is the reciprocity map, and $(U_j)$ is the
filtration in $K^*$ determined by the valuation. A comprehensive
exposition of all these facts is given, e.~g., in \cite{S,
Ch.\,IV, Ch.\,XV}.

However, if $\oK$ is not perfect, there exists no reasonable
theory of upper numbering of ramification subgroups. The ``lower''
ramification subgroups can still be defined, however, the
ramification filtration in the group $G$ does not determine that
in $G/H$. (Examples were given, e.~g., in \cite{L, H}.)

In the present article we treat the class of fields $K$ with $[\oK:\oK^p]=p$. (In particular,
this holds for a two-dimensional local field $K$.)
In the case $\chr K=p$, we work with a relative situation $K/k$,
when a complete subfield $k$ in $K$ with a perfect residue field is supposed
to be fixed. (In the mixed characteristic case, i.~e., when $\chr K=0$, a subfield
$k$ can be chosen in a canonical way.)
 For a Galois extension $L/K$ we
introduce a new lower filtration on $\Gal(L/K)$ indexed by a
special linearly ordered set $\BI$ (see \S1). Then a
Hasse-Herbrand function $\Psi_{L/K}\:\BI\to\BI$ can be defined
with all the usual properties. Therefore, a theory of upper
ramification groups, as well as the ramification theory of
infinite extensions, can be developed.

If we consider abelian extensions of exponent $p$, the
ramification filtration determines a dual filtration on the
additive (resp.  multiplicative) group of $K$ via Artin-Schreier
(resp. Kummer) duality. In the case $e_{K/k}=1$, this dual
filtration is described explicitly in \S2.

In \S3, we consider a field $K$ with a discrete valuation of rank
two. (Main examples are provided by  2-dimensional local or
local-global fields.) We introduce a new index set
$\BI_2\supset\BI$ in this case. This yields  finer lower and upper
filtrations on Galois groups of extensions of $K$.

In the remaining part of the article we deal with an equal
characteristic two-dimensional local field $K$. We introduce
some filtration on the group $K_2\tp K$, which is other than the
filtration induced by the valuation. Our filtration is indexed
by $\BI_2$, and it has a better behavior with respect to the
norm map than the usual filtration.  Finally, we prove that the
reciprocity map of two-dimensional local class field theory (see
\cite{P, F1}) identifies this filtration with the ramification
filtration of \S3.

The ramification filtration constructed in this paper for 2-dimensional
local fields has been generalized to $n$-dimensional local fields
by V.~A.~Abrashkin (see \cite{A1, A2}). He announced a theorem which
is a local version of Grothendieck anabelian conjecture and which
is stated in terms of this filtration.

The author expresses his gratitude to the University of Nottingham
for hospitality, very friendly and stimulating atmosphere.
I am deeply grateful to Prof. I.~B.~Fesenko for many valuable discussions
and for his generous help during the author's stay in Nottingham.

The final version of this paper has been written during the author's
stay at IHES (Bures-sur-Yvette) in 2001. I would like to thank
this institute for its remarkable hospitatlity.

\head
\S0. Definitions, notation and preliminary facts
\endhead

\subhead 0.1. General notation\endsubhead

The letter $p$ always denotes a prime number. This is the
characteristic of residue fields of all discrete valuation
fields under consideration; $v_p(a)$ is the $p$-adic
exponent of a rational or $p$-adic number $a$.

Let $K$ be a field with a discrete valuation of rank 1. Then

$\bullet$ $v_K$ denotes the (normalized) valuation of $K$ as well
as its prolongation  onto the algebraic closure of $K$ (which is
unique provided that $K$ is complete);

$\bullet$ $\CO_K$ is the valuation ring;

$\bullet$ $\gM_K$ is the maximal ideal of $\CO_K$;

$\bullet$ $\oK=\CO_K/\gM_K$;

$\bullet$ $U_K=\CO_K^*$;

$\bullet$ $U_K(m)=\{u\in\CO_K\,|\,v_K(u-1)\ge m\}$, $m\ge1$;

$\bullet$ $\CR_K$ consists of Teichm\"uller representatives of elements
of the maximal perfect subfield in $\oK$;

A finite extension $L/K$ is said to be {\it ferociously ramified}
if $[L:K]=[\oL:\oK]_\insep$.

The group $\BQ^2=\BQ\times\BQ$ is linearly ordered as follows:
$$
(a_1,a_2)<(b_1,b_2)\iff \cases &\text{either } a_2<b_2 \\
                   &\text{or } a_2=b_2 \text{ and } a_1<b_1;
            \endcases
$$

$\bullet$ $\BQ_+=\{i\in\BQ|i>0\}$;

$\bullet$ $\BQ^2_+=\BQ\times\BQ_+$.

\subhead 0.2. Elimination of wild ramification\endsubhead
(See \cite{E, H, Z1, KZ}.)

Let $K$ be a complete discrete valuation field of any
characteristic with the residue field $\oK$ of characteristic
$p>0$.

If $K$ is of characteristic 0, denote by $k$ the set of all
$x\in K$ which are algebraic over the fraction field $k_0$ of
$W(F)$, where $F=\bigcap\oK^{p^i}$. Obviously, $k$ is a complete
subfield of $K$ with perfect residue field, and it is maximal
with respect to this  property.

Next, if $\chr K=p$, we fix a {\it base subfield} $B$ in $K$,
which is complete with respect to the valuation of $K$ and has
$\BF_p$ as a residue field. (The possible base subfields are
exactly all the $\BF_p((\tau))$, where $v_K(\tau)>0$.) In this
case we denote by $k_0$ the completion of $B(\CR_K)$, and by $k$ the algebraic closure
of $k_0$ in $K$.

In both cases $k$ is said to be the {\it constant subfield} of
$K$. We denote by $v_0$ the valuation $K^*\to\BQ$ which
is equivalent to $v_K$ and such that $v_0\mid_{k_0}=v_{k_0}$.
Denote $e_K=e_{K/k_0}=v_0(\pi)^{-1}$, where $\pi$ is any prime of
$K$.

An extension $L/K$ is said to be {\it constant} if $L=lK$ where
$l$ is a certain algebraic extension of $k$. Obviously, in this
case $l$ is the constant subfield of $L$. Notice that finite
separable constant extensions of $K$ are exactly $K(a)/K$ where
$a$ is algebraic and separable over $k$. An extension is said to be {\it almost
constant} if it lies in a compositum of a constant extension and an unramified one.
Equivalently, $L/K$ is almost constant if $K'L/K'$ is a constant extension
for some unramified extension  $K'/K$. Notice the following properties.

1. The compositum of two almost constant extensions is almost constant.
Therefore, one can consider the maximal almost constant subextension in a given
finite extension.

2. Any intermediate extension in an almost constant extension $L/K$ is also almost constant.
Indeed, let $K\subset K'\subset L'\subset L$. We have $L\subset CU$, where
$C/K$ is a constant extension, $U/K$ is an unramified one. Then $L'\subset (CK')(UK')$,
$CK'/K'$ is a constant extension, $UK'/K'$  is an unramified one.

3. If $L/K$ and $M/L$ are almost constant, then $M/K$ is almost
constant.

Indeed, $LU/U$ and $MU'/U'$ are constant for some unramified
extensions $U/K$ and $U'/L$. Let $U''$ be the inertia subfield
in $U'/K$. Then $U'=U''L$, whence $MU'=MU''$. Now $UU''/K$ is
unramified,
$LUU''/UU''$ and $MUU''/LUU''$ are constant, whence $M/K$ is
almost constant.

4. Any tamely ramified extension is almost constant.

We give the name {\it infernal} to the opposite type of
extensions. Namely, a finite extension $L/K$ is said to be
infernal if the only almost constant subextension in $L/K$ is $K$
itself. The following  two properties are proved in the subsection
{\bf1.2} of \cite{KZ}.

1. If $L/K$ is infernal, then any intermediate extension in $L/K$ is infernal.

2. If $L/K$ is infernal, then $K'L/K'$ is also infernal for any
almost constant extension $K'/K$.

We say that $K$ is {\it standard} if $e_{K/k}=1$. It is obvious
that any unramified or constant  or ferociously ramified finite
extension of a standard field is a standard field as well.

Further, $K$ is said to be {\it almost standard} if there exists
an unramified  extension $L/K$ such that $L$ is standard. It is
easy to prove that an almost constant or ferociously ramified
finite extension of an almost standard field is almost standard.

The following statement is a particular case of Theorem 3.2.1 in
\cite{KZ}. However, a simpler proof in this case can be given.

\proclaim{Proposition 0.2.1} Let $K$ be a complete discrete
valuation field, $\chr K=p$, $k$ the constant subfield of $K$,
$L/K$  an infernal  extension. Then there exists a finite purely
inseparable extension $l/k$ such that $lK$ is almost standard, and
$lL/lK$ is ferociously ramified.
\endproclaim

\demo{Proof} Notice that $K$ is a finite extension of a standard
field $K_0$. For example, one can write $K=H((\pi))$ and take
$K_0=H((\pi_0))$, where $\pi_0$ is a prime in the base subfield.
Let $K_1/K_0$ be the maximal almost constant subextension in
$K/K_0$. Then $K/K_1$ is infernal, and $K_1$ is almost standard.
It is sufficient to prove both assertions for $K_2L/K_2K$, where
$K_2/K_1$ is any given unramified extension. Therefore, we may
assume without loss of generality that $K_1$ is standard.

Now it is sufficient to check that $lL/lK_1$ is ferociously
ramified for a certain purely inseparable $l/k_1$, where $k_1$ is
the constant subfield in $K_1$. Using induction on $[L:K_1]$, we
reduce this assertion to the case $[L:K_1]=p$.

First, let $L/K_1$ be inseparable; then $L=K_1(x)$, $x^p=a\in
K_1$. The fact that $L/K_1$ is infernal implies that
$\{v_{K_1}(a-b^p-c)|b\in K_1,c\in k_1\}$ is bounded from above.
Then one can write $a=a_0+b^p+c$ with the maximal possible
$v_{K_1}(a_0)<\infty$. Let $l=k_1(\pi_l)$, where
$\pi_l^p$ is a prime in $k_1$. Then $lL=lK_1(\root p\of{a_0})$
is ferociously ramified over $lK_1$.

Next, assume that $L/K_1$ is
separable and normal. Then $L=K_1(x)$, $x^p-x=a\in K_1$. Let $\pi$
be a prime in $k_1$. Then we may assume $K_1=H_1((\pi))$.

Write $a=a_0+\wp(b)+c$ with $b\in K_1$, $c\in k_1$ with maximal
possible $v(a_0)$. Suppose $v(a_0)\ge0$. Then $L\subset
K_1(x_0,x_c)$, $x_0^p-x_0=a_0$, $x_c^p-x_c=c$, whence $L/K_1$ is
almost constant, a contradiction. Therefore, $v(a_0)=-n<0$.

Let
$s=s_{L/K_1}=\min\{r|\overline{\pi^na_0}\in\bigl(\oK_1\bigr)^{p^r}\}$.
We have $s<\infty$ by the maximality of $v(a_0)$. Now we use
induction on $s$.

If $s_{L/K_1}=0$, then already $L(\pi')/K_1(\pi')$ is ferociously
ramified, where $(\pi')^p=\pi$. If $s\ge1$, it is clear that
$(p,n)=1$. Denote $L'=L(\pi'), K'=K_1(\pi')$, $(\pi')^p=\pi$. We
distinguish two cases.

1. $a_0$ can be written in the form $$ a_0=\sum_{i=-n}^{[-n/p]}
b_i^p\pi^i+a_1, $$ where $b_i\in H$, $v_K(a_1)>-n/p$. Then
$a_0\equiv b_m(\pi')^m+\dots \mod \wp(K')$, where $m=-n/p$, dots
denote terms of higher order, and we obtain
$s_{L'/K'}=s_{L/K_1}-1$.

2. $a_0$ cannot be written in this form. Then we see immediately
that $s_{L'/K'}=0$.

In both cases one applies the assumption of induction to $L'/K'$.

It remains to consider the case when $L/K_1$ is of degree $p$,
separable, but not normal. Let $L_1/K_1$ be the normal closure of
$L/K_1$, $K_2/K_1$ the maximal tamely ramified subextension in
$L_1/K_1$. Then $[L_1:K_2]=p$, $L_1/K_2$ is normal, and $K_2$ is
almost standard. We can apply the already considered case to
$UL_1/U$, where $U$ is a suitable unramified extension of $K_2$.
\qed
\enddemo

\proclaim{Proposition 0.2.2} {\rm (Epp \cite{E})} Let $L/K$ be a
finite extension of complete discrete valuation fields, $k$ the
constant subfield of $K$. Then there exists a finite extension
$l/k$ such that $e_{lL/lK}=1$.
\endproclaim

\demo{Proof} Let $\chr K=p$. Let $K'/K$ be the maximal almost
constant subextension in $L/K$, $k'$ the constant subfield in
$K'$. Then $K'$ is almost standard, and $L/K'$ is infernal. There
exists $l_1/k$ such that $l_1K'/l_1K$ is unramified. By
Proposition 0.2.1, there exists $l_2/k'$ such that $l_2L/l_2K'$ is
ferociously ramified. Let $l=l_1l_2$. Then $lK'/lK$ is unramified,
$lL/lK'$ is ferociously ramified, whence $e(lL/lK)=1$.

In the case $\chr K=0$ the argument is essentially the same.
Instead of purely inseparable extensions $l/k$ one can take
cyclic extensions with sufficiently big ramification
jumps. For a detailed proof, see \cite{Z1, KZ}. \qed
\enddemo

\proclaim{Corollary}
Let $K$ be a complete discrete valuation field, $k$ the constant
subfield of $K$. Then there exists a finite extension $l/k$ such that
$lK$ is standard.
\endproclaim

\subhead 0.3. Two-dimensional local fields\endsubhead
(See \cite{HLF, MZ1, MZ2}.)

We shall freely use the terminology and notation from \cite{HLF, MZ1}.
In particular, a two-dimensional local field is a
field $K$ which is complete with respect to a discrete valuation $v$, and
such that the residue field $\oK$ of $K$ is a complete
discretely valued field with perfect residue field. Throughout
this article we assume that $\oK$ is of characteristic $p$.

A system of local parameters in $K$ is any $(t,\pi)$, where $v(\pi)=1$,
$v(t)=0$, and the residue class of $t$ in $\oK$ is a prime element
of the latter field. The choice of such a system determines a rank
two discrete valuation on $K$, and valuations obtained this way
are equivalent.  The group of principal units of $K$ with
respect to any of these rank two valuations is denoted by $V_K$.
Next, the group of principal units with respect to rank one
valuation of $K$ is denoted by $U_K(1)$.
In other words, $V_K$ (resp. $U_K(1)$) consists of all $a\in\CO_K$
such that the residue class of
$a$ in $\oK$ belongs to $\CO^*_{\oK}$ (resp. equals $1$).

\proclaim{Proposition 0.3}
Let $K$ be a two-dimensional local field (of characteristic
either 0 or $p$), $k$ the constant subfield of $K$. Then $K$
is standard if and only if $K\isom k\<t\>$.
\endproclaim

\demo{Proof}
Consider the embedding $k\<T\>\to K$, which maps $f(T)$ to $f(t)$.
It is natural to denote its image by $k\<t\>$. Then $K/k\<t\>$
is a totally ramified extension with the ramification index $e_{K/k}$.
\qed
\enddemo

\bigbreak
\subhead 0.4. Topological $K$-groups\endsubhead (See \cite{HLF, F3, Z2}.)

\nopagebreak
\proclaim{Well-known lemma 0.4.1}
Let $F$ be a field, $\al,\be\in F^*$, $\al,\be\ne0,1$,
$\al\be\ne\al-1$. Then
$$
\{1-\al,1-\be\}=\bigl\{\al(1-\be),1-(1-\al^{-1})^{-1}
\be\bigr\}
$$
in $K_2F$.
\endproclaim

\demo{Proof}
Using the Steinberg relation, we obtain
$$
\split
\{1-\al,1-\be\}
& = \{1-\al+\al\be,1-\be\}-\Bigl\{1+\frac{\al\be}{1-\al},1-\be\Bigr\}
\\
& = -\{1-\al+\al\be,\al\}-\Bigl\{1+\frac{\al\be}{1-\al},1-\be\Bigr\}
\\
& = -\Bigl\{1+\frac{\al\be}{1-\al},\al\Bigr\}
-\Bigl\{1+\frac{\al\be}{1-\al},1-\be\Bigr\}
\\
& = \Bigl\{\al(1-\be),1+\frac{\al\be}{1-\al}\Bigr\}.
\qed
\endsplit
$$
\enddemo

We recall the definition and some properties of topological
$K$-groups $K_2\tp K$.

Let $K$ be a two-dimensional local field.
The topology on $K_2 K$ is defined as the finest topology such that
the natural map $K^*\times K^*\to K_2 K$ is sequentially
continuous, and subtraction in $K_2 K$ is sequentially continuous.
Then the intersection of all neighborhoods of zero in $K_2K$
coincides with $\cap_{l\ge1}lK_2K$. The group $K_2K/\cap_{l\ge1}lK_2K$
with the quotient topology is denoted by $K_2\tp K$. When it does
not lead to a confusion, the class  of $\{a,b\}\in K_2K$ in $K_2\tp K$ is also denoted
by $\{a,b\}$.

Denote by $\VK_2\tp K$ (resp. $U(1)K_2\tp K$) the subgroup of
$K_2\tp K$ generated by all symbols $\{u,x\}$, where $x\in K$ and
$u\in V_K$ (resp. $u\in U_K(1)$).


Next, following \cite{HLF, Ch. VI} and \cite{F3}, we give an alternative description
of topology on $\VK_2\tp K$.

Let $\tau$ be a topology on a set $X$. The sequential saturation
of $\tau$ is defined as the finest topology $\tau'$ on $X$ such that
the set of convergent sequences for $\tau'$ is the same as that
for $\tau$. For topological spaces $X$, $Y$ the $*$-product
topology on $X\times Y$ is the sequential saturation of the
product topology.

\proclaim{Theorem 0.4.2 (\cite{HLF, 6.6; F3, Theorem 4.6})}
 Let $(t,\pi)$ be local parameters of $K$.
Then the homomorphism
$$
\align
g:V_K\times V_K&\to \VK_2\tp K
\\
(\al,\be)&\mapsto \{\al,\pi\}+\{\be,t\}
\endalign
$$
induces a homeomorphism between $(V_K\times V_K)/\Ker g$
and $\VK_2\tp K$. Here $(V_K\times V_K)/\Ker g$ is endowed
with the quotient topology of $*$-product topology.
\endproclaim

\proclaim{Corollary 0.4.3}
Any neighborhood of 0 in $\VK_2\tp K$ contains $p^n\VK_2\tp K$ for some $n$.
\endproclaim

\demo{Proof}
Let $U$ be a neighborhood of 0 in $\VK_2\tp K$.
If the claim is wrong, for any $i\ge1$ there exists $a_i\in p^i\VK_2\tp K$
such that $a_i\notin U$. Let $a_i=p^ib_i$,
$b_i=\{\al_i,\pi\}+\{\be_i,t\}$, $\al_i,\be_i\in V_K$.
By Lemma 1.6 in \cite{Z2}, any neighborhood of 1 in $K$ contains
$V_K^{p^m}$ for some $m$. Therefore, $\al_i^{p^i},\be_i^{p^i}\to
1$, whence $a_i\to 0$ in $\VK_2\tp K$ by Theorem 0.4.2.\qed
\enddemo

Let $A$ be a subgroup in $\VK_2\tp K$. We denote
by $\Cl A$ the intersection of all open subgroups containing $A$.
This notation is justified by the fact that any closed subgroup
in $\VK_2\tp K$ coincides with the intersection of all open
subgroups containing it, as explained in \cite{F3, Proof of Theorem 4.6}.
Any set of generators of $A$ is said to be a set of
{\it topological generators} for $\Cl A$.
We say also that $A$ is {\it dense} in $\Cl A$.

\proclaim{Proposition 0.4.4}
Let $L/K$ be a finite extension, $A$ any subgroup in $\VK_2\tp K$.
Then $N_{L/K}\Cl A=\Cl N_{L/K}A$.
\endproclaim

\demo{Proof}
By \cite{HLF, 6.8.2}, the norm of any closed subgroup is closed,
and, therefore, $N_{L/K}\Cl A\supset \Cl N_{L/K}A$.

Conversely, let $U'$ be an open subgroup in $\VK_2\tp L$ containing $N_{L/K}A$.
By Corollary 0.4.3, $U'$ contains $p^n\VK_2\tp L$ for some $n$.
Then $N_{L/K}^{-1}(U')$ is an open subgroup $\VK_2\tp K$;
this is explained in \cite{Z2, Corollary  4.4} and in \cite{HLF,
6.8.2}. It follows $\Cl A\subset N_{L/K}^{-1}(U')$, and we
conclude that $N_{L/K}\Cl A\subset \Cl N_{L/K}A$.
\qed
\enddemo


\proclaim{Proposition 0.4.5} Let $K$ be a two-dimensional local field with
the last residue field $F$,
$\chr K=p$. Then any $a\in U(1)K_2\tp K$ can be uniquely written in the form
$$
a=\sum\Sb i\ge1, j\in\BZ\\(j,p)=1\\\th\in B\endSb
c_{ij\th}\{\pi,1+\th\pi^it^j\}+
\sum\Sb i\ge1, j\in\BZ\\(i,p)=1\\(j,p)\ne1\\\th\in B\endSb
c_{ij\th}\{t,1+\th\pi^it^j\},
$$
where $B$ is a certain fixed $\BF_p$-basis of $F$,
$c_{ij\th}\in\BZ_p$.
Here, for any positive integer $n$,
 the set $\{(j,i)|v_p(c_{ij\th})<n\text{ for some }\th\}$ is
admissible, i.~e., $p^n|c_{ij\th}$ for 
all sufficiently small $j<j(i)$ and
any $\th$.
\endproclaim

For the existence and uniqueness of expansion, see
 \cite{P, \S2, Prop. 1 and 3} when $F$ is finite,
and \cite{F2, Prop. 2.4} in the general case.
The remaining statement follows from \cite{MZ2}; we shall not
use it.

\head \S1. General ramification theory\endhead

\subhead 1.1. Definitions \endsubhead

For the whole paper, we denote
$$
\BI=\{-1,0\}\cup\{(\c,i)|i\in\BQ,\,i>0\}\cup\{(\c,\infty)\}
\cup\{(\cf,i)|i\in\BQ,\,i>0\}\cup\{\infty\}.
$$
This is the index set for lower and upper numbering
of ramification subgroups we are going to introduce.
(The letters $\c$ and $\i$ are related to (almost) {\it constant}
and {\it infernal} extensions respectively.)
This set is
linearly ordered as follows:
$$ \gather -1<0<(\c,i)<(\cf,j)<\infty \text{ for any }i,j;
 \\  (\c,i)<(\c,j) \text{ for any }i<j;
\\ (\cf,i)<(\cf,j) \text{ for any }i<j.
\endgather
$$

Let $K$ be a complete discrete valuation field of any
characteristic with the residue field $\oK$ of characteristic
$p>0$. We assume that $[\oK:\oK^p]=p$.

Let $L/K$ be a finite Galois extension, $G=\Gal(L/K)$. For any
$\al\in\BI$ we are going to define a subgroup $G_\al$ in $G$.

We put $G_{-1}=G$, and denote by $G_0$ the inertia subgroup
in $G$, i.~e.,
$$
G_0=\{g\in G|g(a)-a\in\gM_L\text{ for all }a\in\CO_L\}.
$$

Denote by $K_\c/K$ the maximal almost constant
subextension in $L/K$. To introduce subgroups $G_{(\c,i)}=G_{\c,i}$,
we consider first the case when $K_\c/K$ is constant and contains no unramified subextension.
Then $K_\c=lK$, and we have a natural projection
$$
\pr\:\Gal(L/K)\to\Gal(K_\c/K)=\Gal(l/k)=\Gal(l/k)_0,
$$
where $l$ and $k$ are the constant subfields in $L$ and $K$ respectively. Then
we put $G_{\c,i}=\pr^{-1}(\Gal(l/k)_i)$. In the general case take an unramified
extension $K'/K$ such that $K'L/K'$ contains no unramified
subextension, and the maximal almost constant subextension in
$K'L/K'$ (i.~e., $K'K_\c/K'$) is constant.
We put $G_{\c,i}=\Gal(K'L/K')_{\c,i}$. It is easy to see that the choice of $K'$ plays
no role.

Next,
$$
G_{\c,\infty}=\Gal(L/K_\c)=G_{\c,m}.
$$
for $m$ big enough.
Now we turn to the definition of $G_{(\cf,i)}=G_{\i,i}$, $i>0$.

Assume that $K_\c$ is standard and $L/K_\c$ is ferociously
ramified. Let $t\in \CO_L$, $\ot\notin\oL^p$. We define $$
G_{\cf,i}=\{g\in \Gal(L/K_\c)|v_0(g(t)-t)\ge i\} $$ for all $i>0$.

We recall that $v_0=v_{k_0}$; thus, $v_0(L^*)=e_L^{-1}\BZ$.
Notice that in the case under consideration $G_{\cf,i}=G_{ie_L-1}$ in the
notation of \cite{FV}. This follows directly from the definition
of usual ramification groups in \cite{FV, (4.3)}. Therefore,
$G_{\i,i}$ is independent from the choice of $t$.

In the general case choose a finite extension $l'/l$ such that
$l'K_\c$ is standard and $e(l'L/l'K_{\c})=1$, see Proposition
0.2.2 and Corollary. Then it is clear that
$\Gal(l'L/l'K_\c)=\Gal(L/K_\c)$. Next, $l'L/l'K_\c$ is infernal,
whence there is no non-trivial unramified subextension in
$l'L/l'K_\c$. Together with $e(l'L/l'K_\c)=1$ this means that
$l'L/l'K_\c$ is ferociously
ramified. We define
$$G_{\cf,i}=\Gal(l'L/l'K_\c)_{\cf,i}=\Gal(l'L/l'K)_{\cf,i} $$
for all $i>0$.

Finally, put $G_\infty=\{e\}$.

\proclaim{Proposition 1.1.1}
$G_{\cf,i}$ does not depend on the choice of $l'$.
\endproclaim

\demo{Proof}
Let $l'/l$ and $l''/l$ be finite extensions such that
$l'K_\c$ and $l''K_\c$ are standard, and
$$
e(l'L/l'K_{\c})=e(l''L/l''K_{\c})=1.
$$
Then $l'l''K_\c$ is standard, and $e(l'l''L/l'l''K_{\c})=1$.
We may therefore assume that $l'\subset l''$. Let $t\in\CO_{l'L}$,
$\ot\notin(\overline{l'L})^p$. Then we have $t\in\CO_{l''L}$,
$\ot\notin(\overline{l''L})^p$, and we conclude that
$\Gal(l''L/l''K_\c)_{\i,i}=\Gal(l'L/l'K_\c)_{\cf,i}$
immediately by the definition of $G_{\i,i}$ in the ferociously
ramified case.
\qed
\enddemo

Thus, we have defined $G_\al$ for all $\al\in\BI$, and we see immediately
that $\al<\be$ implies $G_\al\supseteq G_\be$. We say that $\al\in\BI$
is a jump for $L/K$ if $G_\al\ne G_\be$ for any $\be>\al$.

\proclaim{Proposition 1.1.2}
For any $\al\in\BI$, $G_\al$ is a normal subgroup in $G$.
\endproclaim

\demo{Proof}
We have nothing to prove, if $\al=-1$ or $\al=\infty$. It is well
known that $G_0$ is normal, see, e.~g., \cite{FV, (4.3)}. Next,
let $\al=(\c,i)$. We may assume that $K'/K$ in the definition of
$G_{\c,i}$ is normal. Let $g\in G$, $h\in G_\al$.
Extend $g$ and $h$ onto $K'L/K$ so that $h\in\Gal(K'L/K')$,
this can be done since $h\in G_\al\subset G_0$.
We see
$$
ghg^{-1}|_{K'}=g|_{K'}\cdot g^{-1}|_{K'}=e,
$$
whence $ghg^{-1}\in\Gal(K'L/K')=G_0$.

Let $K'_\c/K$ be the maximal almost constant subextension in
$K'L/K'$, and let $k'$, $l'$ be the constant subfields in $K'$ and
$K'L$ respectively. We have $g(k')=k'$, whence $g(K_\c')=K_\c'$.
Now
$$
\pr(ghg^{-1})=g|_{K_\c'}\cdot h|_{K_\c'}\cdot g^{-1}|_{K_\c'}\in
g\Gal(l'/k')_ig^{-1}=\Gal(g(l')/g(k'))_i=\Gal(l'/k')_i.
$$
This proves also the case $\al=(\c,\infty)$.

In the case $\al=(\i,i)$ the same argument is applicable. Indeed,
$l'/l$
in the definition may be chosen to be normal, and the relation
$v_0(g(t)-t))\ge i$ is stable with respect to any automorphisms of
discretely valued fields.
\qed
\enddemo

\subhead 1.2. Compatibility with subgroups\endsubhead
\proclaim{Proposition 1.2}
Let $H$ be a subgroup in $G$. Then for any $\al\in\BI$ we have
$H_\al=H\cap G_\al$.
\endproclaim

\demo{Proof}
The assertion is obvious for $\al=-1,0,\infty$. Next, denote by $M$
the intermediate field in $L/K$ which corresponds to $H$.
Then $K_\c M/M$ is almost constant, whence $K_\c M\subset M_\c$.
On the other hand, $L/K_\c M$ is infernal, and we conclude $K_\c M= M_\c$.
This proves the assertion for $\al=(\c,\infty)$.
Next, if $e(\tl L/\tl K_\c)=1$, then $e(\tl L/\tl M_\c)=1$ as well.
It is clear now that
$$
H_{\cf,i}=\Gal(L/M_\c)\cap G_{\cf,i}=H_{\c,\infty}\cap G_{\cf,i}=
H\cap G_{\cf,i}
$$
for all $i>0$.

Finally, let $\al=(\c,i)$. Denote by $K'/K$ any unramified
extension such that $K'L/K'$  contains no unramified
subextension, and the maximal almost constant subextension in
$K'L/K'$ is constant.
Then $K'L/K'M$ possesses these two properties as well.
Let $k'$, $l'$, $m'$ be constant subfields in $K'$, $K'L$, $K'M$
respectively. Then the natural projection
$\Gal(K'L/K'M)\to\Gal(l'/m')$ is the restriction of the natural
projection $\pr:\Gal(K'L/K')\to\Gal(l'/k')$.
 Therefore,
$$
\align
H_\al&=H\cap\pr^{-1}(\Gal(l'/m')_i)
\\
&=H\cap\pr^{-1}(\Gal(l'/m')\cap\Gal(l'/k')_i)
\\
&=H\cap\Gal(K'L/m'K')\cap\Gal(K'L/K')_\al
\\
&=H\cap G_\al,
\endalign
$$
the latter equality follows from
$$
H\cap\Gal(K'L/K')_\al\subset
H\cap\Gal(K'L/K')=\Gal(K'L/K'M)\subset\Gal(K'L/m'K').
\qed
$$
\enddemo

\subhead 1.3. Hasse-Herbrand functions and compatibility with
quotient groups\endsubhead

To compute the ramification of a quotient group $G/H$, we have to
express the ramification number of given $\si\in G/H$ in terms of
ramification numbers of its representatives in $G$. For constant
extensions, we use \cite{Se, Ch.~IV, Prop.~3 and Lemma 5}. Now we
carry out the same calculation for the case of ferociously
ramified $L/K$.

Let $t\in\CO_L$, $\ot\notin\oL^p$. Let
$$
i_G(g)=v_0(g(t)-t)=\max\{i|g\in G_{\cf,i}\}
$$
for $g\in G=\Gal(L/K)$, $g\ne1$.

\proclaim{Lemma 1.3.1} Let $L/K$ be ferociously ramified, $H$ a
normal subgroup in $G$, $\si\in G/H$, $\si\ne1$. Then $$
i_{G/H}(\si)=i_G(g_1)+\dots+i_G(g_n), $$ where $g_1,\dots,g_n$ are
all the representatives of $\si$ in $G$.
\endproclaim

\demo{Proof}
Similarly to \cite{Se, Ch.~IV, Prop.~3}.
\qed
\enddemo

Let $L/K$ be an arbitrary finite Galois extension, $K_\c/K$ the maximal
almost constant subextension in $L/K$.
We define $\Phi_{L/K}\:\BI\to\BI$ as follows:
$$
\split
\Phi_{L/K}(\al)&=\al  \text{ for } \al=-1,0,(\c,\infty),\infty;
 \\
\Phi_{L/K}((\c,i))&=
 \biggl(\c,\frac1{e(L/K)}\int_0^i|\Gal(K_\c/K)_{\c,t}|dt\biggr)
\text{ for all } i>0;
 \\
\Phi_{L/K}((\cf,i))&=
\biggl(\cf,\int_0^i|\Gal(L/K)_{\cf,t}|dt\biggr)
\text{ for all } i>0.
\endsplit
$$
It is easily seen that $\Phi_{L/K}$ is bijective and increasing;
we introduce also $\Psi_{L/K}=\Phi_{L/K}^{-1}$.

\proclaim{Proposition 1.3.2}
Let $L/K$ and $M/K$ be finite Galois extensions with $M\subset L$.
Then $\Phi_{L/K}=\Phi_{M/K}\circ\Phi_{L/M}$.
\endproclaim

\demo{Proof} The proposition follows immediately from the definitions
in the following three cases: $L/K$ is almost constant; $L/K$ is
infernal (by Lemma 1.3.1); $M/K$ is almost constant whereas $L/M$
is infernal.

Let now $M/K$ be infernal, and $L/M$ be almost constant. Then
there exists an almost constant Galois extension $L'/K$ such that
$L=L'M$. (The argument is the same as in the proof of the 1st
property of infernal extension in {\bf0.2}.) The extension $L'/K$
is in fact the maximal almost constant subextension in $L/K$. By
the definition $\Phi_{M/K}=\Phi_{L/L'}$. It is also easy to see
that $\Phi_{L'/K}=\Phi_{L/M}$. These two functions commute since
$\Phi_{M/K}$ affects only the ``infernal'' part of $\BI$, and
$\Phi_{L'/K}$ only the ``constant'' one. Therefore, $$
\Phi_{M/K}\circ\Phi_{L/M}=\Phi_{L/L'}\circ
\Phi_{L'/K}=\Phi_{L'/K}\circ\Phi_{L/L'}=\Phi_{L/K}. $$

In the general case let $T/M$ be the maximal
almost constant subextension in $L/M$, and $S/K$ that in $M/K$.
This is clear that $T/S$ is a Galois extension. Let $T'/S$ be the maximal
almost constant subextension in it. Then $T'$  is obviously the maximal
almost constant subextension also in $L/K$ and, therefore, normal over $K$.
We have now
$$
\split
\Phi_{L/K}&=\Phi_{L/T'}\circ\Phi_{T'/K}\\
&=\Phi_{L/T}\circ\Phi_{T/T'}\circ\Phi_{T'/S}\circ\Phi_{S/K}\\
&=\Phi_{L/T}\circ\Phi_{T/S}\circ\Phi_{S/K}\\
&=\Phi_{L/T}\circ\Phi_{T/M}\circ\Phi_{M/S}\circ\Phi_{S/K}\\
&=\Phi_{L/M}\circ\Phi_{M/K}.\qed
\endsplit
$$
\enddemo

\proclaim{Proposition 1.3.3}
Let $G=\Gal(L/K)$, $H$ a normal subgroup in $G$, $M$ the
corresponding intermediate field. Then
for any $\al\in\BI$ we have $(G/H)_\al=G_{\Psi_{L/M}(\al)}H/H$.
\endproclaim

\demo{Proof}
The cases $\al=-1$, $\al=\infty$ are trivial.

Let $\al=0$. If $g\in G_0$, then $g$ acts trivially on
$\oM\subset\oL$, whence $gH\in(G/H)_0$. Conversely, let
$gH\in(G/H)_0$. Then the image of $g$ in $\Gal(\oL/\oK)$ is
$\oh\in\Gal(\oL/\oM)$. If $h$ is a preimage of $\oh$ in
$H=\Gal(L/M)$, then $gh^{-1}\in G_0$, whence $gH\in G_0H/H$.

Let $\al=(\c,i)$, $0<i\le\infty$. Consider first the case when
$K_\c/K$ is constant and contains no unramified subextensions. Let
$k,l,m$ be the constant subfields of $K,L,M$ respectively.
Consider the commutative diagram of natural epimorphisms
$$
\CD
\Gal(L/K) @>\pr_L>> \Gal(l/k)
\\
@VQVV @VqVV
\\
\Gal(M/K) @>\pr_M>> \Gal(m/k)
\endCD
$$
It is sufficient to prove that
$Q^{-1}((G/H)_\al)=G_{\Psi_{L/M}(\al)}H$. We have
$$
\split
Q^{-1}((G/H)_\al) &= Q^{-1}(\pr_M^{-1}(\Gal(m/k)_i))
\\&=\pr_L^{-1}(q^{-1}(\Gal(m/k)_i))
\\&=\pr_L^{-1}(q^{-1}(\Gal(l/k)_{\psi_{l/m}(i)}\Gal(l/m)/\Gal(l/m)))
\\&=\pr_L^{-1}(\Gal(l/k)_{\psi_{l/m}(i)}\Gal(l/m))
\\&=G_{\c,\psi_{l/m}(i)}\cdot\Gal(L/mK)
\\&=G_{\Psi_{L/M}(\al)}\cdot\Gal(L/mK)
\\&\supset G_{\Psi_{L/M}(\al)}H
\endsplit
$$

It remains to prove that $\Gal(L/mK)\subset
G_{\Psi_{L/M}(\al)}H$. Let $g\in \Gal(L/mK)$, then
$\pr_L(g)\in\Gal(l/m)$. Let $h$ be any preimage of $\pr_L(g)$ in
$H=\Gal(L/M)$; then $gh^{-1}\in\Gal(L/lK)=G_{\c,\infty}\subset
G_{\psi_{L/M}(\al)}$.

In the general case take $K'/K$ as in the proof of Prop. 1.2;
$H'=\Gal(K'L/K'M)$. Then $\Psi_{L/M}=\Psi_{K'L/K'M}$, and
$$
\Gal(K'M/K')_\al=\Gal(K'L/K')_{\Psi_{L/M}(\al)}H'/H'.
$$
The natural embedding $\Gal(K'M/K')\hookrightarrow G/H$ identifies
these groups with $(G/H)_\al$ and $G_{\Psi_{L/M}(\al)}H/H$
respectively.

Finally, let $\al=(\i,i)$. Denote by $K_\c^L$ and $K_\c^M$ the
maximal almost constant extensions in $L/K$ and $M/K$
respectively; we have $K_\c^M\subset K_\c^L$. Consider the case
when both $K_\c^L$ and $K_\c^M$ are standard, and $L/K_\c^M$ is
ferociously ramified. In this case we use Lemma 1.3.1 and argue
exactly as in \cite{Se, Ch. IV, Lemmas 3,4,5}. In the general case
take a finite extension $l'/l$ such that $l'K_\c^M$ and $l'K_\c^L$
are standard, and $e(l'L/l'K_\c^M)=1$. Then
$G_{\Psi_{L/M}(\al)}=\Gal(l'L/l'K)_{\Psi_{L/M}(\al)}$, and
$(G/H)_\al=\Gal(l'M/l'K)_\al$. Since
$\Psi_{L/M}(\al)=\Psi_{l'L/l'M}(\al)$ for $\al\ge(\c,\infty)$, we
have $\Gal(l'M/l'K)_\al=\Gal(l'L/l'K)_{\Psi_{L/M}(\al)}H'/H'$,
where $H'=\Gal(l'L/l'M)$; in $G/H$ this yields
$(G/H)_\al=G_{\psi_{L/M}(\al)}H/H$
\qed
\enddemo

\subhead 1.4. Upper numbering\endsubhead

As usual, we denote $G^\al=G_{\Psi_{L/K}(\al)}$ for any $\al\in\BI$.
Then Proposition 1.3.3 can be reformulated as follows.

\proclaim{Corollary 1.4} Let $G$ be the Galois group of a finite
Galois extension of $K$, $H$ a normal subgroup in $G$. Then
for any $\al\in\BI$ we have $(G/H)^\al=G^\al H/H$.\qed
\endproclaim

Let $L/K$ be an arbitrary Galois extension, $L=\bigcup L_i$, $L_i/K$
finite Galois extensions, $G=\Gal(L/K)$. Then we define
$$
G^\al=\lim_\leftarrow\Gal(L_i/K)^\al;\quad L^\al=L^{G^\al}.
$$

Further, one can define $K_\al=(K^{\sep})^\al$ for all $\al\in\BI$.
Then $K_{-1}=K$, $K_0=K^{ur}$, $K_{\c,1}=K^{tr}$, $K_{\c,\infty}=K^{ur}k^\sep$,
$K_\infty=K^\sep$, where $k$ is the constant subfield in $K$.

\subhead1.5. Examples\endsubhead

In all examples $K=F((\pi))$ is of characteristic $p$, standard,
with $\pi\in k$; $t$ is an element of $F\setminus F^p$. Assume $v_0(\pi)=1$.

1. Let $L$ be the splitting field of $X^p-X-\pi^{-pi}t$, where
$i>0$.  Then the only jump of lower filtration for $L/K$ is
$(\i,i)$ and the only jump of upper filtration is $(\i,pi)$.

Indeed, $L/K$ is ferociously ramified, and $\CO_L=\CO_K[\pi^ix]$,
where $x$ is a root of $X^p-X-\pi^{-pi}t$.

2. Let $L$ be the splitting field of $X^p-X-\pi^{-i}t$ where
$i>0$, and $i$ is prime to $p$. Then the only jump of
lower filtration is $(\i,i/p)$ and the only jump of upper
filtration is $(\i,i)$.

In this case $L/K$ is infernal but totally ramified. To compute
the ramification, consider $l'=k(\pi_1)$, $\pi_1^p=\pi$.
Then $l'L/l'K$ is ferociously ramified, and the minimal polynomial
of this extension is $X^p-X-\pi_1^{-pi}t$. Now we apply the
previous example and take into account that $v_0(\pi_1)=p^{-1}$.

3. Let $L$ be the splitting field of $X^p-X-\pi^{-i}t^{p^j}$
where $i>0$, and $i$ is prime to $p$. Then the only
jump of lower filtration is $(\i,i/p^{j+1})$ and the only jump of
upper filtration is $(\i,i/p^j)$.

In this case take $l'=k(\pi_{j+1})$, $\pi_{j+1}^{p^{j+1}}=\pi$. Then
$\pi^{-i}t^{p^j}\equiv\pi_{j+1}^{-pi}t\mod\wp(l'K)$. Therefore,
$l'L=l'K(x)$, where $x$ is a root of $X^p-X-\pi_{j+1}^{-pi}t$.
We conclude that $l'L/l'K$ is ferociously ramified, and the jumps
can be computed as in the first example.

4. Let $L$ be the splitting field of $X^p-X-\pi^{-i}$
where $i>0$, and $i$ is prime to $p$. Then the only
jump of lower filtration is $(\c,i)$ and the only jump of
upper filtration is $(\c,i)$.

\head \S2. Example: abelian extensions of exponent $p$\endhead

We start
with the case $\chr K=p$.
In this case the group of characters of $\Gal(K)$ of exponent
$p$ can be identified via Artin-Schreier theory with the group
$K/\wp(K)$, where $\Gal(K)$ is the absolute Galois group of $K$,
and $\wp\:K\to K$ is the Artin-Schreier map
$x\mapsto x^p-x$. Since the group $\Gal(K)$ possesses a
decreasing filtration indexed by $\BI$, the additive group
of $K$ acquires an increasing filtration indexed by $\BI$.

Let $\al\in\BI$. Denote by $A_\al=A_{\al,K}$ the set of such
$a\in K$ that the polynomial $X^p-X-a$ completely splits in
$K_\al$.  These $A_\al$ form an increasing filtration on $K$:
$$
\wp(K)=A_{-1}\subseteq
A_0\subseteq A_{\c,1}\subseteq\dots\subseteq A_\infty=K.
$$
It
is obvious that $A_0=\CO_K+\wp(K)$.  Further, it is easy to
obtain that $A_{\c,i}=\gM_k^{-i}+A_0$ for all integral $i>0$,
and $A_{\c,\infty}=k+A_0$.
As for $A_{\i,i}$, these subgroups can be calculated
explicitly only in case $K$ is standard.

\proclaim{Proposition}
Let $K$ be standard, i.~e., $K=F((\pi))$ with $\pi\in k$. Then
$$
A_{\i,i}=k+A_0+\sum_{l=0}^\infty \pi^{-p^li}F^{p^l}[[\pi]].\tag1
$$
\endproclaim

\demo{Proof}
Let $a$ belong to the right hand side of (1), $x$
be a root of $X^p-X-a$. To prove that $K(x)\subset K_{\i,i}$,
it is sufficient to consider the following cases.

1. $a\in k$, then $K(x)/K$ is constant, whence $K(x)\subset
K_{\c,\infty}\subset K_{\i,i}$.

2. $a\in \CO_K$, then $K(x)/K$ is unramified, whence
$K(x)\subset K_0\subset K_{\i,i}$.

3. $a\in \pi^{-p^li}F^{p^l}[[\pi]]$. In the field $K'=K(\pi')$,
where ${\pi'}^{p^l}=\pi$, we have $a=b^{p^l}$, $b\in
{\pi'}^{-p^li}F[[\pi']]$. Observe that $b^{p^l}\equiv
b\mod\wp(K')$, whence $K'(x)=K'(x')$, ${x'}^p-x'=b$. Clearly,
$v_0(x')=p^{-1}v_0(b)\ge-p^{-1}i$. Without loss of generality we
may assume that $K(x)/K$ is infernal. By Proposition 0.2.1, there exists
 a finite constant
extension $K''/K'$ such that $K''(x)/K''$ is ferociously ramified.
If $g$ is a generator of $\Gal(K''(x)/K'')$, we have $$
v_0(g(x')-x')=0, $$ whence the only jump in the lower filtration
is $\le(\i,p^{-1}i)$, and the only jump in the upper filtration is
$\le(\i,i)$. We conclude $K(x)\subset K_{\i,i}$, and $a\in
A_{\i,i}$.

Conversely, suppose that $a\in A_{\i,i}$, $x^p-x=a$. If
$K(x)/K$ is almost constant, then $a=a_u+a_c$,
such that $K(x_u)/K$ is unramified and $K(x_c)/K$ is constant,
$x_c^p-x_c=a_c\in K$, $x_u^p-x_u=a_u\in K$. It follows
$a_c\in k+\wp(K)$, $a_u\in \CO_K+\wp(K)$, and we are done.

It remains to consider the case of infernal $K(x)/K$.
Let
$$a=b_m\pi^{-m}+b_{m-1}\pi^{-(m-1)}+\dots,$$
where
$b_j\in F$.  Adding an element of $k\subset
A_{\c,\infty}\subset A_{\i,i}$ to $a$, we may assume
$b_j=d_j^{p^{l_j}}$, where $d_j\notin F^p$.
Similarly, adding an element of $\wp(K)$, we may assume
that for any $j$ either $p\nmid j$ or $l_j=0$.
Consider quotients $\varkappa_j=e_{K}^{-1}j/p^{l_j}$.
($e_{K}$ plays the role of
absolute ramification index.) Due to the previous remark, the
values of $\varkappa_j$ are all different. ($\varkappa_r=\varkappa_s$
would imply that at least one of the fractions $r/p^{l_r}$
and $s/p^{l_s}$ is reducible.)
Notice that $\varkappa=\max(\varkappa_j)$ remains invariant
when we replace $K$ with $K(\root p\of\pi)$.

Now let $K'=K(\pi')$, where ${\pi'}^{p^N}=\pi$, $N=\max(l_j)+1$.
Then it is easy to see that $K'(x)=K'(x')$, $$
{x'}^p-x'={\pi'}^{-\varkappa e_{K'}}d+\dots, $$ $d\notin F^p$.
Since $p\mid\varkappa e_{K'}$, we conclude that $K'(x)/K$ is
ferociously ramified, and the jump of ramification in upper
numbering is $\varkappa$. This means $\varkappa\le i$, and $a$
lies in the right hand side of (1). \qed
\enddemo

Now we turn to the case where $\chr K=0$ and $K$ contains all
$p$th roots of unity.

Let $\al\in\BI$. Denote by $B_\al$ the set of such $b\in K^*$
that $\root p\of b\in K_\al$. These $B_\al$ form an increasing
filtration in $K^*$:
$$
(K^*)^p=B_{-1}\subseteq
B_0\subseteq B_{\c,1}\subseteq\dots\subseteq B_\infty=K^*.
$$
It is obvious that $B_0=U_{\frac{pe}{p-1},K}\cdot(K^*)^p$, where
$e=e_K$.  Further, it is easy to obtain that 
$$ 
B_{\c,i}=\cases
U_{\frac{pe}{p-1}-i,k}B_0, &i<\frac{pe}{p-1}, \\ 
k^*B_0,
&i\ge\frac{pe}{p-1}. \endcases 
$$

If $K$ is standard then $B_{\cf,i}$ is generated by $k^*B_0$ and all elements
of the kind $1+\pi^ja^{p^m}$ where $j\ge0$, $a\in\CO_K$ and $m$
is the minimal non-negative integer such
that
$$ \frac{pe}{p-1}-j\le p^{m}i.  $$

\head \S3. Fine ramification theory\endhead

Assume that the algebraic closure of the residue field of $K$ is also endowed
with a valuation $w\:(\oK^\sep)^*\to\BQ$ such that 
the restriction of $w$ to $\oK$ is discrete. It is clear that 
the restriction of $w$ to the maximal perfect subfield of $\oK$ is trivial. 
The main example is given by a two-dimensional
local field $K$. In this case $\oK$ is the field of Laurent power series over a perfect field $F$.
One can take for $w$ the extension of the valuation on $\oK$ onto $\oK^\sep$. Another
example is a two-dimensional local-global field, i.~e., $\oK$ is a field of algebraic functions
in one variable over a perfect field. (Then there is a lot of non-equivalent valuations $w$.)

Now we define the canonical (up to the choice of a base subfield and of $w$)
rank 2  valuation $\bv_0$ on $K$.

\proclaim{Proposition 3.1}
Let $k$ be the constant subfield in $K$.
There exists a unique valuation $\bv_0=(v_1,v_2)\:(K^\alg)^*\to\BQ\times\BQ$ such that
$v_2(a)=v_0(a)$ for all $a\in K^*$, $v_1(u)=w(\ou)$ for all $u\in U_{K^\alg}$,
and $\bv_0(c)=(0,v_0(c))$ for all $c\in k$.
The value group of $\bv_0|_{K^*}$ is isomorphic to $\BZ\times \BZ$.
\endproclaim

\demo{Proof}
Let $\pi_0$ be a prime element of $k_0$. For $a\in K^\alg$, let
$v_0(a)=\frac mn$, $m,n\in\BZ$. Then $\pi_0^{-m}a^n\in
U_{K^\alg}$, and we see that the only possible value for $\bv_0(a)$
is $\bigl(w\bigl(\overline{\pi_0^{-m}a^n}\bigr)/n,m/n\bigr)$. 
The value group of $\bv_0|_{K^*}$ is
generated by $\bv_0(\pi)$ and $(1,0)$, where $\pi$ is a prime of
$K$.
\qed
\enddemo

Introduce the index set
$$
\BI_2=\BI\cup\{(i_1,i_2)|i_1,i_2\in\BQ,i_2>0\}.
$$
We extend the ordering of $\BI$ onto $\BI_2$ assuming
$$
(\cf,i_2)<(i_1,i_2)<(i_1',i_2)<(\cf,i_2')
$$
for all $i_2<i_2'$, $i_1<i_1'$.

We turn to the definition of $G_{i_1,i_2}$ where $G$ is the Galois
group of a given finite Galois extension $L/K$. Assume first that
$K_\c$ is standard and $L/K_\c$ is ferociously ramified. Let $t\in
\CO_L$, $\ot\notin\oL^p$. We define $$ G_{i_1,i_2}=\biggl\{g\in
\Gal(L/K_\c)\bigg| \bv_0\biggl(\frac{g(t)}t-1\biggr)\ge(i_1,i_2)\biggr\} $$
for $i_1,i_2\in\BQ$, $i_2>0$. Observe that $\CO_L=\CO_{K_\c}[t]$,
whence $\bv_0(g(t)t^{-1}-1)\ge(i_1,i_2)$ implies $\bv_0(g(a)a^{-1}-1)\ge(i_1,i_2)$
for any $a\in\CO_L$, $a\ne0$. This proves the independence of $G_{i_1,i_2}$
from the choice of $t$.

In the general case we choose $l'/l$
such that $l'K_\c$ is standard and $l'L/l'K_\c$ is ferociously
ramified and put $$ G_{i_1,i_2}=\Gal(l'L/l'K_\c)_{i_1,i_2}. $$
Like in \S1, we see that $G_{i_1,i_2}$ does not depend on the
choice of $l'$.
(Therefore, only $k_0$ and $w$ are involved.)

We say that $\al\in\BI_2$ is a jump for $L/K$ if $G_\al\ne G_\be$ for
any $\be>\al$. It is easy to see that any jump is either $-1$, or $0$,
or $(\c,i)$ with integral $i$, or $(i_1,i_2)$.

Like in {\bf 1.3}, one constructs Hasse-Herbrand functions
$\Phi_{2,L/K}\:\BI_2\to\BI_2$ and $\Psi_{2,L/K}=\Phi_{2,L/K}^{-1}$ which extend
$\Phi$ and $\Psi$ respectively. Namely,
$$
\align
\Phi_{2,L/K}((i_1,i_2))&=\int_{(0,0)}^{(i_1,i_2)}|\Gal(L/K)_{t}|dt
\\
&=\sum_{i=1}^{m}(h_{i}-h_{i-1})|\Gal(L/K)_{h_i}|+((i_1,i_2)-h_m)|\Gal(L/K)_{i_1,i_2}|,
\endalign
$$
where $h_1<\dots<h_m$ are all jumps for $L/K$ between $(0,0)$ and
$(i_1,i_2)$, and $h_0=(0,0)$.

The assertions of {\bf 1.2}, {\bf 1.3}, {\bf 1.4} all remain valid, if one replaces $\BI$,
$\Phi$, $\Psi$ with $\BI_2$, $\Phi_2$, $\Psi_2$ respectively. In what follows,
we shall write $\Phi$ and $\Psi$ instead of $\Phi_2$ and $\Psi_2$.

\head \S4. The subgroups $S_\al$ in $U(1)K_2\tp K$ \endhead

In the rest of this paper $K$ denotes a two-dimensional local
field (see {\bf0.3}) such that
$\chr K=p$.  We fix a base subfield in $K$.  The constant
subfield of $K$ is denoted by $k$.

We fix a diskrete valuation of renk one on the first residue field of $K$. It has a unique extension
to the valuation $w\:\bigl(\oK^\alg\bigr)^*\to\BQ$. We do not require that $w(\oK^*)\subset\BZ$.
Then a valuation $\bv_0=(v_1,v_2)\: (K^\alg)^*\to\BQ^2$ is defined, see \S3.
Let $t,\pi$ be local parameters of $K$. Introduce a matrix
$$
\bde_{\pi,t}=\pmatrix v_1(t) & v_2(t) \\ v_1(\pi) & v_2(\pi)\endpmatrix^{-1}.
$$
Then it is easy to see that the valuation $\bv_K\:K^*\to\BZ^2$, which is associated with $\pi$ and $t$, can be written as $\bv_K=\bv_0|_K\cdot \bde_{\pi,t}$. In most cases, we shall not
mention explicitly that we have fixed certain local parameters $t,\pi$ in $K$.

Sometimes we shall use the notation $U_{\al,K}$ for $\{u\in K|\bv_0(u-1)\ge\al\}$.

Introduce subgroups
$$
\split
& Q_\al=Q_{\al,K}=
\langle\{c,u\}\,|\,c\in k,u\in K, \bv_0(u-1)\ge\al\rangle
\subset U(1)K_2\tp K, \\
& Q_\al^{(n)}=\{a\in U(1)K_2\tp K \,|\, p^na\in Q_\al\}, \\
& S_\al=S_{\al, K}=\Cl\bigcup_{n\ge0} Q_{p^n\al}^{(n)},
\endsplit
$$
where $\al\in\BQ^2_+$, $n\ge0$.

To check that $S_\al$ is a subgroup in $U(1)K_2\tp K$, observe that
 $pQ_\al\subset Q_{p\al}$.
Therefore, $Q_\al^{(n)}\subset Q_{p\al}^{(n+1)}$ for all $n$.


\proclaim{Proposition 4.1}
Let $K$ be standard, $\pi, t$ local parameters of $K$ such that $\pi$ is a constant.
Then the following elements of $U(1)K_2\tp K$ topologically generate $S_\al$:
$$
\aligned
& p^{r_{ij}}\{\pi,1+\th\pi^it^j\},\quad \th\in B, \,(j,p)=1;\\
& p^{r_{ij}}\{t,1+\th\pi^it^{j}\},\quad 
\th\in B, \,(i,p)=1,\,(j,p)\ne1,\,j\ne0,
\endaligned
\tag1
$$
where $B$ is a certain fixed $\BF_p$-basis of the last residue field of $K$,
$r_{ij}=r_{ij}(\al)$ is the minimal non-negative integer 
such that $p^{r_{ij}}(j,i)\ge p^{v_p(j)}\al\bde_{\pi,t}$.
\endproclaim

\comment

\demo{Proof}
The group $Q_\al$ is topologically generated by various
$\{\pi,(1+\th\pi^i t^j)^{p^r}\}$, where $p^r(j,i)\ge\al\bde_{\pi,t}$,
and either $(i,p)=1$, or $(j,p)=1$.
We may assume here $j\ne0$ (since $K_2\tp k=0$, see \cite{Kn, Th\'eor\`eme 1}).
Let $(i,p)=1$. Then
$$
\split
\{\pi,1+\th\pi^it^j\}&=i^{-1}\{\pi^i,1+\th\pi^it^j\}
\\
&=i^{-1}\{-\th t^j,1+\th\pi^it^j\}
\\
&=i^{-1}j\{t,1+\th\pi^it^j\}
\endsplit.
$$
We conclude that
$Q_\al$ consists of all elements $a\in U(1)K_2\tp K$ such that in the standard expansion
of Proposition 0.4.5
$$
a=\sum\Sb i\ge1, j\in\BZ\\(j,p)=1\\\th\in B\endSb
c_{ij\th}\{\pi,1+\th\pi^it^j\}+
\sum\Sb i\ge1, j\in\BZ\\(i,p)=1\\(j,p)\ne1\\\th\in B\endSb
c_{ij\th}\{t,1+\th\pi^it^j\}
$$
we have $v_p(c_{ij\th})\ge v_p(j)+r_{ij}$.
\qed
\enddemo

\endcomment

\proclaim{Corollary 1}
In the setting of Proposition 4.1, the subgroup $pS_{\al/p}$,
together with the following elements of $U(1)K_2\tp K$,
topologically generate $S_\al$:
$$
\split
& \{\pi,1+\th\pi^it^j\},\quad \th\in B, \,p\nmid j,\, (j,i)\ge\al\bde_{\pi,t};\\
& \{t,1+\th\pi^it^{j}\},\quad \th\in B, \,p\nmid i,\,p\mid j,\,(j,i)\ge p^{v_p(j)}\al\bde_{\pi,t}.
\endsplit
$$
\endproclaim

(One can say that these elements are topological generators of
$S_\al$ modulo $pS_{\al/p}$.)

\proclaim{Corollary 2}
In the setting of Proposition 4.1, let
$$
 a_{ij\th}=\cases \{\pi,u_{ij\th}\}, & \th\in B, \,p\nmid j,\, (j,i)\ge\al\bde_{\pi,t};\\
 \{t,u_{ij\th}\}, & \th\in B, \,p\nmid i,\,p\mid j,\,(j,i)\ge p^{v_p(j)}\al\bde_{\pi,t}.
\endcases
$$ where $u_{ij\th}\equiv 1+\th\pi^i t^j\mod U_K(j+1,i)$.
Assume that all $a_{ij\th} \in S_\al$ and that $u_{ij\th}$ form
a part of some system of topological generators of $K^*$. Then
all these $a_{ij\th}$ are topological generators of $S_\al$
modulo $pS_{\al/p}$.
\endproclaim

Next, for a fixed $t$, for any $j\ne0$ and $i>0$ denote
$$
W_{ij}=W_{ij,K}=\Bigl\{1+\sum_{\nu=1}^\infty c_\nu t^{\nu j}
\big| c_\nu\in k,\, v_k(c_\nu)\ge\nu i\Bigr\}.
$$
Obviously, $W_{ij}$ is a subgroup in $K^*$.

\proclaim{Proposition 4.2}
Let $K=k\<t\>$, $\pi$ a prime element in $k$,
$\al\in\BQ^2_+$, $j\ne0$, $i>0$,
$(j,i)\ge p^{v_p(j)}\al\bde_{\pi,t}$.
Then for any $u\in W_{ij}$ we have 
$\{t,u\}\in S_\al$. 
\endproclaim

\demo{Proof of Propositions 4.1 and 4.2}
Denote by $A_\al^0$ the subgroup of $U(1)K_2\tp K$ generated
by all elements \thetag1 and put $A_\al=\Cl A_\al^0$. Denote also 
$B_d=\{\th^{p^d}|\th\in B\}$, $d\ge0$. This is also an $\BF_p$-basis of $F$.

\smallskip
{\it Step 1}: $A_\al\subset S_\al$.

Let $(i,p)=1$, $(j,p)\ne1$. Then
$$
\aligned
\{\pi,1+\th\pi^it^j\}&=i^{-1}\{\pi^i,1+\th\pi^it^j\}
\\
&=i^{-1}\{-\th t^j,1+\th\pi^it^j\}
\\
&=i^{-1}j\{t,1+\th\pi^it^j\},
\endaligned
\tag2
$$
whence $p^{v_p(j)}\{t,1+\th\pi^it^j\}\in Q_\be$,
$\be=(j,i)\bde_{\pi,t}^{-1}$.
If $r_{ij}\le v_p(j)$, we have
$$
p^{r_{ij}}\{t,1+\th\pi^it^j\}\in Q_\be^{(v_p(j)-r_{ij})}
\subset S_{p^{r_{ij}-v_p(j)}\be}\subset S_\al.
$$
If $r_{ij}> v_p(j)$, we have
$$
p^{r_{ij}}\{t,1+\th\pi^it^j\}\in p^{r_{ij}-v_p(j)}Q_\be
\subset Q_{p^{r_{ij}-v_p(j)}\be}\subset S_\al.
$$

Similarly, for $(j,p)=1$ we have
$$
p^{r_{ij}}\{\pi,1+\th\pi^it^j\}\in p^{r_{ij}}Q_\be
\subset Q_{p^{r_{ij}}\be}\subset S_\al.
$$
Therefore, $A_\al^0\subset S_\al$, and $A_\al\subset S_\al$.

\smallskip
{\it Step 2}: the elements
$p^{r_{ij}}\{\pi,1+\th\pi^it^j\}$ and $p^{r_{ij}}\{t,1+\th\pi^it^j\}$
are in $A_\al$ for all $i>0$, all $j\ne0$, and all $\th\in B_d$, 
$d=\min(v_p(i),v_p(j))$.

If $(i,p)=1$, $(j,p)\ne1$, the symbol $p^{r_{ij}}\{\pi,1+\th\pi^it^j\}$
is a multiple of $p^{r_{ij}}\{t,1+\th\pi^it^j\}$,
see \thetag2. The case $(j,p)=1$ is similar.
If  $d=\min(v_p(i),v_p(j))>0$, $i=p^di_0$ and $j=p^dj_0$,
we have $r_{i_0j_0}=r_{ij}$, and for $\tau=\pi,t$:
$$
p^{r_{ij}}\{\tau,1+\th\pi^it^j\}=
p^d\cdot p^{r_{i_0j_0}}\{\tau,1+\th_0\pi^{i_0}t^{j_0}\}\in A_\al.
$$

\smallskip
{\it Step 3}: elements 
$\{t,u\}$ 
(from Prop. 4.2) belong to $A_\al$.

Any element of $W_{ij}$ can be written as an infinite product
of $1+\th\pi^mt^{nj}$ with $m\ge ni>0$ and $\th\in B_d$,
$d=\min(v_p(m),v_p(ni))$.
We have
$$
(nj,m)\ge(nj,ni)\ge np^{v_p(j)}\al\bde_{\pi,t}
\ge p^{v_p(n)}p^{v_p(j)}\al\bde_{\pi,t}
=p^{v_p(nj)}\al\bde_{\pi,t}.
$$
Therefore, $r_{m,nj}=0$, whence 
$\{t,1+\th\pi^mt^{nj}\}$
are in $A_\al$ by Step 2.

\comment

It remains to check that $\{u,u'\}\in S_\al$ for
$u=1+\th\pi^mt^{nj}$, $u'=1+\th'\pi^{m'}t^{n'j}$,
where $\th,\th'\in B$, $m\ge ni$, $m'\ge n'i$.
In view of Lemma 0.4.1, we have
$$
\{u,u'\}=
\{-\th\pi^mt^{nj}(1+\th'\pi^{m'}t^{n'j}),w\},
$$
where
$$
w=1+\frac{\th\th'\pi^{m+m'}t^{(n+n')j}}{1+\th\pi^mt^{nj}}
\in W_{ij}\cap U_K(m+m').
$$
The symbol $\{-\th,w\}$ is infinitely divisible, and, therefore,
equals $0$ in $K_2\tp K$. We conclude that
$$
\{u,u'\}=m\{\pi,w\}+nj\{t,w\}+\{u',w\}.
$$
We already know that $\{\pi,w\}$ and $\{t,w\}$. Applying the
same transformation to $\{u',w\}$, and taking into account
that $v(u'-1)+v(w-1)=m+2m'>v(u-1)+v(u'-1)$, we expand $\{u,u'\}$
into a convergent sum of elements of $A_\al$.

\endcomment

\smallskip
{\it Step 4}: $Q_\al\subset A_\al$ (for all $al$)
implies $S_\al\subset A_\al$.

Indeed, $A_\al$ consists of all $a\in U(1)K_2\tp K$ such that
in the standard expansion of $a$ as in Proposition 0.4.5
we have $v_p(c_{ij\th})\ge r_{ij}(\al)$ whenever $j\ne0$, and
$c_{ij\th}=0$ when $j=0$. Since $\max(r_{ij}(p^n\al),n)=r_{ij}(\al)+n$,
we obtain
$$
A_{p^n\al}\cap p^nU(1)K_2\tp K=p^nA_\al.
$$
Therefore, $Q_\al\subset A_\al$ implies $Q_{p^n\al}^{(n)}\subset A_\al$.

\smallskip
{\it Step 5}: $Q_\al\subset A_\al$.

The group $Q_\al$ is topologically generated by $\{\pi,1+\th\pi^it^j\}$
and $\{1+\eta\pi^r,1+\th\pi^it^j\}$, where $\th,\eta\in B$,
$(j,i)\ge\al\bde_{\pi,t}$, and $r$ is any positive integer.
We have to check that these elements are in $A_\al$.

If $j=0$, then all these elements are equal to $0$ since
$K_2\tp k=0$ (see \cite{Kn, Th\'eor\`eme 1}). Thus, we may assume
$j\ne0$. We may also reduce to the case $p\nmid(i,j)$. Indeed, if
$\min(v_p(i),v_p(j))=d>0$, we have
$$
\{\pi,1+\th\pi^it^j\}
=p^d\{\pi,1+\th^{p^{-d}}\pi^{i_0}t^{j_0}\}
\in p^dA_{\al/p^d}\subset A_\al,
$$
and the same for $\{1+\eta\pi^r,1+\th\pi^it^j\}$.

If $(j,p)=1$, we have $r_{ji}(\al)=0$, and $\{\pi,1+\th\pi^it^j\}$
is just one of the elements \thetag1. If $(j,p)\ne1$, we have
 $r_{ji}(\al)\le v_p(j)$ and, by \thetag2,
$\{\pi,1+\th\pi^it^j\}$ is a multiple of
$p^{r_{ij}}\{t,1+\th\pi^it^j\}$ which is one of the elements \thetag1.

Finally, by Lemma 0.4.1,
$$
\aligned
\{1+\eta\pi^r,1+\th\pi^it^j\}
&=-\{1+\th\pi^it^j,1+\eta\pi^r\}
\\
&=-\{-\th\pi^it^j(1+\eta\pi^r),w\},
\endaligned
\tag3
$$
where
$$
w=1+\frac{\th\eta\pi^{i+r}t^j}{1+\th\pi^it^j}\in W_{ij}\cap U_K(i+r).
$$
The symbol $\{-\th,w\}$ is infinitely divisible, and, therefore,
equals $0$ in $K_2\tp K$. We conclude that
$$
\{1+\eta\pi^r,1+\th\pi^it^j\}=-i\{\pi,w\}-j\{t,w\}-\{1+\eta\pi^r,w\}.
$$
We already know from the previous paragraph that $\{\pi,w\}\in A_\al$.
Next, $j\{t,w\}$ is a multiple of $p^{v_p(j)}\{t,w\}$, and
$\{t,w\}\in A_{\al/p^{v_p(j)}}$ by Step 3. Therefore,
$j\{t,w\}\in p^{v_p(j)}A_{\al/p^{v_p(j)}}\subset A_\al$.
Applying the same transformation \thetag3 to $\{1+\eta\pi^r,w\}$,
we expand $\{1+\eta\pi^r,1+\th\pi^it^j\}$
into a convergent sum of elements of $A_\al$.
\qed
\enddemo

\head \S5. Behavior of $S_\al$ in certain types of
extensions  \endhead

\subhead5.1\endsubhead
In this section we assume that $K$ is an equal characteristic two-dimensional local field
with the following additional property.
\smallskip
(*) {\it Let $L/K$ be  any finite unramified extension. Then the
extension $l/k$ is also unramified, where $k$ and $l$ are the
constant subfields of $K$ and $L$ respectively. }
\smallskip
The class of such fields is stable with respect to constant or
unramified extensions. However, an almost constant extension of
$K$ need not be of this type: adjoin a root of $X^p-X-\pi^{-1}-t$
to $\BF_p((\pi))\<t\>$. Note that any standard field satisfies
(*).

This condition on $K$, together with Proposition 0.2.1,
easily implies that a certain constant purely inseparable extension of $K$ is standard.

The aim of this section is to prove the following three assertions.

\proclaim{Proposition 5.1.1} Let $L/K$ be an unramified
extension, and let $\al\in \BQ^2_+$. Then we have
$N_{L/K}S_{\al,L}=S_{\al,K}$.
\endproclaim

\proclaim{Proposition 5.1.2} Let $L/K$ be a constant totally
ramified extension, and let $\al\in \BQ^2_+$. Then
$N_{L/K}S_{\al,L}=S_{\al,K}$.
\endproclaim


\proclaim{Proposition 5.1.3} Let $K$ be standard, $L/K$ a cyclic
ferociously ramified extension of degree $p$ with the
ramification jump $h$ in lower numbering, and let
$\al\in\BQ^2_+$. Then:

(1) $N_{L/K}S_{\al,L}=S_{\al+(p-1)h,K}$ if $\al>h$;

(2)$N_{L/K}S_{\al,L}$ is a subgroup in $S_{p\al,K}$ if $\al\le h$;

(3) $(S_{p\al,K}:N_{L/K}S_{\al,L})=p$ if $\al\le h$ and the last residue field of $K$
is quasi-finite.
\endproclaim

For the purposes of the proof, it is convenient to introduce notation
$S_{\al,K}^{0}=\bigcup_{n\ge0}Q_{p^n\al}^{(n)}$. Then $S_{\al,K}=\Cl S_{\al,K}^0$,
and to prove, say, Proposition 5.1.1, it is sufficient to check
that its conditions imply that $N_{L/K}S_{\al,L}^0$ is dense in
$S_{\al,K}$.

\subhead 5.2\endsubhead
Let $L/K$ be a finite extension. Since there is no torsion in $\VK_2\tp K$,
the natural map $\VK_2\tp K\to\VK_2\tp L$ is injective, and we may
identify $\VK_2\tp K$ with a subgroup in $\VK_2\tp L$.

\proclaim{Lemma 1}
Let $L/K$ be inseparable of degree $p$. Then $\VK_2\tp K\subset p\VK_2\tp L$.
\endproclaim

\demo{Proof}
Let $L/K$ be totally ramified. Let $\pi,t$ be local parameters of $L$. Then
$\pi^p,t$ are local parameters of $K$. In $\VK_2\tp K$ one can choose
a system of topological generators of type $\{\pi^p,1+\th(\pi^p)^it^j\}$, $p\nmid j$,
and $\{t,1+\th(\pi^p)^it^j\}$, $p\mid j$.
Evidently, all these generators lie in $p\VK_2\tp L$.

For a ferociously ramified $L/K$, one should change the roles
of $\pi$ and $t$. \qed
\enddemo

\proclaim{Lemma 2}
Let $L/K$ be an inseparable constant extension of degree $p$.
Then $S^0_{\al,K}\subset S^0_{p\al, L}.$
\endproclaim

\demo{Proof}
It is sufficient to check that $Q_{\al,K}\subset Q_{p\al,L}$.
Let $l$ and $k$ be the constant subfields in $L$ and $K$ respectively.
Then $k=l^p$,
$$
\split
Q_{\al,K} & = \langle\{c^p,u\}\,|\,c\in l,u\in K, \bv_0(u-1)\ge\al\,\}\rangle\\
&\subset \langle\{c^p,u\}\,|\,c\in l,u\in L, \bv_0(u-1)\ge\al\,\}\rangle\\
&=\langle\{c,u^p\}\,|\,c\in l,u\in L, \bv_0(u-1)\ge\al\,\}\rangle\\
&\subset \langle\{c,u\}\,|\,c\in l,u\in L, \bv_0(u-1)\ge p\al\,\}\rangle\\
&=Q_{p\al,L}.\qed
\endsplit
$$
\enddemo

\proclaim{Proposition 5.2}
Let $L/K$ be a constant purely inseparable totally
ramified extension, and let $\al\in \BQ^2_+$. Then
$N_{L/K}S_{\al,L}=S_{\al,K}$.
\endproclaim

\demo{Proof}
We may assume $[L:K]=p$. Let $\pi$ be a prime in the constant subfield of $L$. We have
$$
\split
N_{L/K}(pQ_{\al,L})
&=N_{L/K}(\langle\{c^p,u\}\,|\,c\in l,u\in L, \bv_0(u-1)\ge\al\,\}\rangle) \\
&=\langle\{c^p,u^p\}\,|\,c\in l,u\in L, \bv_0(u-1)\ge\al\,\}\rangle\\
&\subset \langle\{c^p,u\}\,|\,c\in l,u\in K, \bv_0(u-1)\ge p\al\,\}\rangle \\
&=Q_{p\al,K},
\endsplit
$$
whence $N_{L/K}Q_{\al,L}\subset S_{\al,K}$, and  $N_{L/K}S_{\al,L}\subset S_{\al,K}$
for all $\al\in\BQ^2_+$. Conversely, let $a\in S^0_{\al,K}$. Then $a=pb=N_{L/K}b$,
$b\in\VK_2\tp L$ by Lemma 1. By Lemma 2, $pb\in S_{p\al,L}^0$, whence $b\in S_{\al,L}^0$.
\qed
\enddemo

\subhead 5.3. Proof of Proposition 5.1.1\endsubhead

By Proposition 0.2.1 and condition (*), there exists an inseparable constant
extension $K'/K$ such that $K'$ is standard. In view of Proposition 5.2,
we may assume without loss of generality that $K'=K=k\<t\>$, $L=l\<t\>$.
Fix a prime $\pi\in k$, it is also a prime element of $l$.

To prove that $N_{L/K}S_{\al,L}\subset S_{\al,K}$, it is sufficient
to check that the norms of all standard topological generators of
$N_{L/K}S_{\al,L}$ (as in Prop. 4.1) belong to $S_{\al,K}$. Since
$p^rS_{\al/p^r,K}\subset S_{\al,K}$, we may consider only the case
$r_{ij}=0$. Now, for $(j,i)\ge\al\bde_{\pi,t}$,
$$
N_{L/K}\{\pi,1+\th\pi^it^j\}=\{\pi,N_{L/K}(1+\th\pi^it^j)\}\in Q_{\al,K}
\subset S_{\al,K}
$$
by the definition of $Q_{\al,K}$.
Next, for $(j,i)\ge p^{v_p(j)}\al\bde_{\pi,t}$,
$$
N_{L/K}\{t,1+\th\pi^it^j\}=\{t,N_{L/K}(1+\th\pi^it^j)\}\in S_{\al,K}.
$$
Indeed, all the conjugates of $1+\th\pi^it^j$ are in $W_{ij,L}$,
whence $N_{L/K}(1+\th\pi^it^j)\in W_{ij,L}\cap K=W_{ij,K}$, and
we may apply Prop. 4.2.

Conversely,
to prove that $S_{\al,K}\subset N_{L/K}S_{\al,L}$,
it is sufficient
to check that all standard topological generators of $S_{\al,K}$
belong to $N_{L/K}S_{\al,L}$ (see Prop. 0.4.4). Again, we can easily
reduce to the case $r_{ij}=0$. We have $1+\th\pi^it^j\in N_{L/K}U_L(j,i)$,
whence, for $(j,i)\ge\al\bde_{\pi,t}$,
$$
\{\pi,1+\th\pi^it^j\}\in N_{L/K} Q_{\al,L}\subset N_{L/K} S_{\al,L}.
$$
To see that $\{t,1+\th\pi^it^j\}\in N_{L/K}S_{\al,L}$
 for $(j,i)\ge p^{v_p(j)}\al\bde_{\pi,t}$,
observe that the norms of topological generators of $W_{ij,L}$
topologically generate $W_{ij,K}$, whence $N_{L/K}W_{ij,L}=W_{ij,K}$.
\qed

\subhead 5.4. Proof of Proposition 5.1.2\endsubhead

Let $L/K$ be a constant totally ramified extension. Denote by $\tilde L/K$
the normal closure of $L/K$,
by $K'/K$ the maximal unramified subextension in
$\tilde L/K$ and by $L'/L$ that in $\tilde L/L$.
By Proposition 5.1.1, $S_{\al,L}=N_{L'/L}S_{\al,L'}$,
and $S_{\al,K}=N_{K'/K}S_{\al,K'}$, and it is sufficient
to prove Proposition for $\tilde L/L'$ and $\tilde L/K'$.

Therefore, we may assume that $L/K$ is normal, and, further,
that $L/K$ is cyclic of prime  degree.
As in {\bf 5.3}, there exists an inseparable constant extension $K'/K$
such that $K'$ is standard. In view of Proposition 5.2,
we may assume without loss of generality that $K'=K$.

Let $\pi,t$ be local parameters in $L$ such that $\pi$ is a constant.
Then $\pi_0,t$ are local
parameters of $K$, where $\pi_0=N_{L/K}\pi$; we have $K=k\<t\>$.

We start with the case $[L:K]=p$. Let the ramification jump of $L/K$
be $(\c,h)$.

First, we prove that $N_{L/K}S_{\al,L}\subset S_{\al,K}$.
Notice that the group
$\{\{\pi,u\}|\bv_0(u-1)\ge\al\}$ is topologically generated by
$$ \split
&\{\pi,1+\th\pi_0^it^j\},\quad (j,pi)\ge\al\bde_{\pi,t}, \\
&\{\pi,1+\th\pi^it^j\},\quad (i,p)=1,\,(j,i)\ge\al\bde_{\pi,t}.
\endsplit
$$
It follows from Prop. 4.1 that $S_{\al,L}$ is topologically generated by
these elements and all $p^{r_{ij}}\{t,1+\th\pi^it^j\}$ of Prop. 4.1.
As in {\bf 5.3}, it is sufficient to consider only those $(j,i)$ where
$r_{ij}=0$. Notice also that $(i,p)=1$ implies that
$\{\pi,1+\th\pi^it^j\}$ is a multiple of $\{t,1+\th\pi^it^j\}$
(see \thetag2 in \S4).

Therefore, $N_{L/K}Q_{\al,L}$ is topologically generated by
$$
N_{L/K}\{\pi,1+\th\pi_0^it^j\}=\{\pi_0,1+\th\pi_0^it^j\},\quad (j,pi)\ge\al\bde_{\pi,t} \text{ (i.~e., $(j,i)\ge\al\bde_{\pi_0,t}$)},
$$
and
$$
\split
N_{L/K}\{t,1+\th\pi^it^j\}
&=\{t,N_{L/K}(1+\th\pi^it^j)\} \\
&=\Bigl\{t,1+\sum_{\nu=1}^p c_\nu t^{\nu j}\Bigr\},
\quad (i,p)=1,\,(j,i)\ge v_p(j)\al\bde_{\pi,t}.\endsplit
$$
Here $c_\nu\in k$, $v_k(c_\nu)\ge\frac1p(\nu i+(p-1)h)$,
$\nu=1,\dots,p-1$; $v_k(c_p)=i$.
Then Proposition 4.2 implies that
$N_{L/K}\{t,1+\th\pi^it^j\}\in S_{\al,K}$ (note that
 $(j,i/p)\ge v_p(j)\al\bde_{\pi_0,t}$).

It remains to prove that $N_{L/K}S_{\al,L}$ is dense in $S_{\al,K}$.
In view of Proposition 0.4.3, it is sufficient
to show that $N_{L/K}S_{\al,L}+pS_{\al/p,K}$ is dense in $S_{\al,K}$.
To do this, we shall construct a system of topological generators
$a_{ij\th}$ for $S_{\al,K}/pS_{\al/p,K}$ as in Corollary 2 to Proposition 4.1
such that $a_{ij\th}\in N_{L/K}S_{\al,L}$.
In this argument we shall also make additional
requirement for all $a_{ij\th}$ with $v_p(i)=0$.
Let $j=p^mj_0$, $m=v_p(j)$.
We require that $a_{ij\th}=\{t,u_{ij\th}\}$, where
$$
u_{ij\th}=1+\sum_{\nu=p^m}^\infty c_\nu t^{\nu j_0}
\eqno{(1)}
$$
with $c_\nu\in k$, $v_k(c_\nu)\ge \frac i{p^m}\nu$.

Apply induction on $v_p(j)$. If $v_p(j)=0$, take $a_{ij\th}=N_{L/K}\{\pi,1+\th\pi_0^it^j\}$ for all $(j,i)\ge\al\bde_{\pi_0,t}$.

Suppose that $a_{ij\th}$ with $v_p(j)\le n$ have already been constructed.
Take a pair $(j,i)$ with $v_p(j)=n+1$
and $(j,i)\ge p^{n+1}\al\bde_{\pi_0,t}$.
Consider
$$
\split
N_{L/K}\{t,1+\th\pi^it^{j/p}\}&=\{t,N_{L/K}(1+\th\pi^it^{j/p})\} \\
&=
\{t,1+\sum_{r=1}^{p-1}c_rt^{rj/p}+\th^p\pi_0^it^{j}\},
\endsplit
$$
where $c_r\in k$, $v_k(c_r)\ge ri+(p-1)h$.
We take $a_{ij\eta^p}=N_{L/K}\{t,1+\eta\pi^it^{j/p}\}$ if $(j,i)<(0,h)$.
In the remaining case it is easy to see that for some integers
$m_{rs\th}$, $r=1,\dots,p-1$, $s\ge v_k(c_r)$, $\th\in B$
we have
$$
N_{L/K}(1+\eta\pi^it^{j/p})\prod_{r,s,\th}u_{rj/p,s,\th}^{m_{rs\th}}
\equiv 1+\eta^p\pi_0^it^j \mod U_{j+1,i,K}.
$$
The left-hand side of this relation can be then taken as $u_{ij\eta^p}$. It is seen immediately
that it satisfies (1), and we put $a_{ij\eta^p}=\{t,u_{ij\eta^p}\}$.

Then all these $a_{ij\th}$ belong to $N_{L/K}S_{\al,L}$, and they topologically generate
the quotient group $S_{\al,K}/pS_{\al/p,K}$ by Corollary 2 to
 Proposition 4.1.

It remains to consider the case $[L:K]=q\ne p$. Since $Q_{\al,K}=qQ_{\al,K}$,
we have
$$
S_{\al,K}=qS_{\al,K}=N_{L/K}S_{\al,K}\subset N_{L/K}S_{\al,L}.\qed
$$
The converse can be proved as in {\bf5.3}.
\qed

\subhead 5.5. Proof of Proposition 5.1.3\endsubhead

Let $\pi$ and $t$ be local parameters in $L$.
Then $\pi$ and $t_0=N_{L/K}t$ are local parameters in $K$.

First, we prove that $N_{L/K}S_{\al,L}\subset S_{p\al,K}$ for $\al\le h$.
It is sufficient to check that
$N_{L/K}(S_{\al,L}/pS_{\al/p,L})\subset S_{p\al,K}/pS_{\al,K}$.
The quotient group $S_{\al,L}/pS_{\al/p,L}$ is generated
by the classes of $N_{L/K}\{\pi,1+\th\pi^i t^j\}$,
$(j,i)\ge\al \bde_{\pi,t}$, and $N_{L/K}\{t,1+\th\pi^i t_0^j\}$,
$(j,i)\ge p^{v_p(j)+1}\al \bde_{\pi,t}$, use Proposition 4.2. We have
$$
\split
N_{L/K}\{\pi,1+\th\pi^i t^j\}&=\{\pi,N_{L/K}(1+\th\pi^i t^j)\}\in S_{p\al,K}, \\
N_{L/K}\{t,1+\th\pi^i t_0^j\}&=\{t_0,1+\th\pi^i t_0^j\}\in S_{p\al,K}.
\endsplit
$$
Similarly one shows $N_{L/K}S_{\al,L}\subset S_{\al+(p-1)h}$ for $\al\ge h$.

Next, we have $N_{L/K}Q_{\al,L}=Q_{\al+(p-1)h,K}$ for all $\al>h$.
Let $n>0$. Denote $K_n=k\<t_0^{p^n}\>\subset K$,
$L_n=k\<t^{p^n}\>\subset L$, where $k$ is the constant subfield of
$K$. Then $L_n/K_n$ is a ferociously ramified extension with the
ramification jump $p^nh$. It follows $$
p^nN_{L/K}Q_{p^n\al,L}^{(n)}=N_{L_n/K_n}Q_{p^n\al,L_n}=Q_{p^n\al+(p-1)p^nh,K_n}
=p^nQ_{p^n\al+p^n(p-1)h,K}^{(n)}, $$
 whence
$N_{L/K}Q_{p^n\al,L}^{(n)}=Q_{p^n\al+p^n(p-1)h,K}^{(n)}$. Since all $Q_{p^n\al,L}^{(n)}$
topologically generate $S_{\al,K}$, we conclude $N_{L/K}S_{\al,L}=S_{\al+(p-1)h,K}$ for all $\al>h$.

It remains to prove the third part of Proposition.
Assume that the last residue field
of $K$ is quasi-finite. Then the map
$$
N_{L/K}\:U_L(j,i)/U_L(j+1,i)\to U_K(j,pi)/U_K(j+1,pi)
$$
has the cokernel of order $p$, where $(j,i)=\al\bde_{\pi,t}$.
Let $u\in U_K(j,pi)$ generate this cokernel.
Then $\{\pi,u\}\notin N_{L/K}K_2\tp L$,
and $\{\pi,u\}\in Q_{h,K}\subset S_{\al,K}$ for any $\al\le h$.
This shows $(S_{p\al,K}:N_{L/K}S_{\al,L})\ge p$ for $\al\le h$.

Application  of subfields $K_n$ and $L_n$ yields an easy calculation
of all $N_{L/K}\{\pi,1+\th\pi^i t^j\}$ and $N_{L/K}\{t,1+\th\pi^i t^j\}$. We obtain that
$N_{L/K}S_{\al,L}$ together with $\{\pi,u\}$ generates $S_{p\al,K}$ for all $\al\le h$,
whence $(S_{p\al,K}:N_{L/K}S_{\al,L})\le p$ for $\al\le h$.
\qed

\head \S6. A filtration on $K_2\tp K$ and reciprocity map.
 \endhead

In this section $K$ is a two-dimensional local field with a quasi-finite
residue field, $\chr K=\chr\oK=p$.

\subhead6.1\endsubhead
For any $\al\in\BI_2$ we introduce a subgroup $Fil_\al K_2\tp K\subset K_2\tp K$
so that $Fil_\be K_2\tp K\subset Fil_\al K_2\tp K$ whenever $\al<\be$. It is easy to see that
for some unramified extension  $\tilde K/K$ the field $\tilde K$ satisfies the condition
(*), see \S5. Denote

$Fil_\infty K_2\tp K=0$;

$Fil_\al K_2\tp K=N_{\tilde K/K}S_{\al,\tilde K}$ for $\al\in\BQ^2_+$;

$Fil_{\i,\al_2} K_2\tp K=\Cl\bigcup\limits_{\al_1\in\BQ}Fil_{\al_1,\al_2} K_2\tp K$ for $\al_2\in\BQ_+$;

$Fil_{\c,\infty} K_2\tp K=T_K:=\Cl\bigcup\limits_{\al\in\BQ^2_+}Fil_{\al}K_2\tp K$;

$Fil_{\c,i}K_2\tp K=T_K+\{\,\{t,u\}\,|\,u\in k,\, v(u-1)\ge i\}$ for all $i\in\BQ_+$
 if $K=k\<t\>$ is standard;

$Fil_{\c,i}K_2\tp K=N_{K'/K}Fil_{\c,i}K_2\tp K'$, if $K$ has the property (*)
and $K'$ is a purely inseparable constant
extension such that $K'$ is standard;

$Fil_{\c,i}K_2\tp K=N_{\tilde K/K}Fil_{\c,i}K_2\tp \tilde K$ in
the general case, where $\tilde K/K$ is an unramified extension
such that $\tilde K$ satisfies (*);


$Fil_0K_2\tp K=U(1)K_2\tp K+\{\,\{t,\th\}\,|\,\th\in\CR_K\}$;

$Fil_{-1}K_2\tp K=K_2\tp K$.

\proclaim{Proposition 6.1.1}
$Fil_\al K_2\tp K$ are well defined.
\endproclaim

\demo{Proof} Independence from
the choice of $\tilde K$ follows from Proposition 5.1.1 and from
 $N_{k''/k'}U_{i,k''}=U_{i,k'}$ in an unramified $k''/k'$.  If
$u'\in U_K$, $u\in U_k(1)$, then the definition of $Q_\al$ implies
$\{u',u\}\in T_K$. This shows
independence from the choice of $t$. Independence from the
choice of $K'$ follows from two observations:

(1) $N_{K''/K'}T_{K''}=T_{K'}$ if $K''=k''\<t\>$, $K'=k'\<t\>$,  $k''/k'$ is purely inseparable;

(2) in this case $N_{k''/k'}U_{i,k''}=U_{i,k'}$.\qed
\enddemo

\proclaim{Proposition 6.1.2}
$\bigcup_{\al>0}Fil_{\al}K_2\tp K=U(1)K_2\tp K$.
\endproclaim

\demo{Proof}
This is clear if $K$ is standard. Next,
$N_{K'/K}U(1)K_2\tp K'=U(1)K_2\tp K$ if $K'/K$ is a purely inseparable or unramified extension.
\qed
\enddemo

\remark{Remark}
If $K$ is standard, then obviously
$$
Fil_{\c,i}K_2\tp K=T_K+U(i)K_2\tp K
$$
for a positive integer $i$.
\endremark

\subhead6.2. Norm map in a purely inseparable constant
extension\endsubhead

\proclaim{Proposition 6.2}
Let $L/K$ be a finite purely inseparable constant extension.
Then $N_{L/K}Fil_\al K_2\tp L=Fil_\al K_2\tp K$ for all $\al\in\BI_2$.
\endproclaim

\demo{Proof}
For $\al=(\i,\al_2)$,
or $\al\in\BQ^2_+$, see Proposition 5.1.2. The cases $\al=-1,
(\c,\infty),\allowmathbreak\infty$ are trivial. For
$\al=(\c,i)$, it is sufficient to apply the observations
in the proof of Proposition 6.1.1. Finally, let
$\al=0$. We have $$ N_{L/K}U(1)K_2\tp L=U(1)K_2\tp K, $$ and $$
N_{L/K}(\{t,\th\})=\{t,\th^{p^n}\}
$$
for all $\th\in\CR_K=\CR_L$, where $p^n=[L:K]$. Therefore,
$$
N_{L/K}Fil_0 K_2\tp L=Fil_0 K_2\tp K.\qed
$$
\enddemo

\subhead6.3. Norm map in a cyclic extension of prime
degree\endsubhead

\proclaim{Proposition 6.3}
Let $L/K$ be a cyclic extension of prime degree, $G=\Gal(L/K)$, $h\in\BI_2$
the only jump of ramification for $L/K$, i.~e., $G_h=G$. Then:

1. If $\al>h$, then $N_{L/K}Fil_\al K_2\tp L=Fil_{\Phi_{L/K}(\al)}K_2\tp K$.

2. If $\al\le h$, then $N_{L/K}Fil_\al K_2\tp L$ is a subgroup
in $Fil_{\Phi_{L/K}(\al)}K_2\tp K$ of index  $[L:K]$.
\endproclaim

\demo{Proof}

1. Reduction to the case when both $K$ and $L$ are almost standard. Choose purely inseparable
$k'/k$ such that $K'=k'K$ and $L'=k'L$ are almost standard.  It is easy to see that $\Phi_{L'/K'}=\Phi_{L/K}$.
Suppose Proposition is valid for $L'/K'$.

Let $\al>h$. Then
$$
\split
N_{L/K}Fil_\al K_2\tp L&=N_{L/K}N_{L'/L}Fil_\al K_2\tp L'\\
&=N_{K'/K}Fil_{\Phi_{L'/K'}(\al)}K_2\tp K'\\
&=N_{K'/K}Fil_{\Phi_{L/K}(\al)}K_2\tp K'\\
&=Fil_{\Phi_{L/K}(\al)}K_2\tp K.
\endsplit
$$

Let $\al\le h$. Then $$ Fil_{\Phi_{L'/K'}(\al)}K_2\tp K'\supset
N_{L'/K'}Fil_\al K_2\tp L' \supset Fil_\be K_2\tp K' $$ for any
$\be>\Phi_{L'/K'}(\al)$. Applying $N_{K'/K}$ to the above
formula, we obtain in view of Proposition 6.2: $$
Fil_{\Phi_{L/K}(\al)}K_2\tp K\supset N_{L/K}Fil_\al K_2\tp
L\supset Fil_\be K_2\tp K $$ for any $\be>\Phi_{L/K}(\al)$. By
class field theory, $K_2\tp K'/N_{L'/K'}K_2\tp L'$ is cyclic of
prime degree. Therefore, $K_2\tp K'$ is generated by
$N_{L'/K'}K_2\tp L'$ and $a$, where $a\notin Fil_\be K_2\tp K'$
for any $\be>\Phi_{L/K}(\al)$. Therefore, $K_2\tp K$ is
generated by $N_{K'/K}a$ and $N_{K'/K}N_{L'/K'}K_2\tp L'$. Since
$$ N_{K'/K}a\in N_{K'/K}Fil_{\Phi_{L/K}(\al)}K_2\tp
K'=Fil_{\Phi_{L/K}(\al)}K_2\tp K', $$ we conclude that
$(Fil_{\Phi_{L/K}(\al)}K_2\tp K:N_{L/K}Fil_\al K_2\tp L)=[L:K]$.

\smallskip

2. The case, when $K$ is almost standard, and $h=-1$, i.~e., $L/K$ is unramified.
Then for $\al\ge(\c,\infty)$ the assertion follows directly from the definition
of $Fil_\al$ with use of Proposition 0.4.4. The case of $\al=(\c,i)$
also follows immediately from the definition of $Fil_\al$.
Further,
$$
\split
N_{L/K}Fil_0 K_2\tp L &= N_{L/K}(U(1)K_2\tp L+\{t\}\cdot\CR_L) \\
&= N_{L/K}U(1)K_2\tp L+\{t\}\cdot N_{L/K}\CR_L)\\
&= U(1)K_2\tp K+\{t\}\cdot \CR_K.
\endsplit
$$

Finally, $N_{L/K}Fil_{-1}K_2\tp L=N_{L/K}K_2\tp L$ is a subgroup of index $[L:K]$
in $K_2\tp K$ by class field theory.

The above argument also reduces Proposition to the case when both $K$ and $L$
are standard. In fact, since $K$ and $L$ are almost standard, there exists an
unramified extension
$K'/K$ such that $K'$ and $K'L$ are standard.
We can therefore identify $Fil_\al K_2\tp K$ with $S_{\al,K}$,
and $Fil_\al K_2\tp L$ with $S_{\al,L}$.

\smallskip

3. Let $K$ be standard, $h=(\c,H)$, $H>0$, or $h=H=0$.
In both cases $L/K$ is a constant extension. Let $k$ and $l$ be the
constant subfields of $K$ and $L$ respectively.
Then $l/k$ is a cyclic extension with the jump of ramification $H$.
By Proposition 4.1.2, we obtain
$N_{L/K}Fil_\al K_2\tp L=Fil_\al K_2\tp K$ for all $\al\ge(\c,\infty)$.
Let $\al=(\c,i)$ or $\al=0=i$. Then
$$
N_{L/K}Fil_\al K_2\tp L = N_{L/K}T_L + N_{L/K}(\{t\}\cdot U_{i,l})=
 T_K+\{t\}\cdot N_{L/K}U_{i,l}.
$$
It remains to notice that $N_{L/K}U_{i,l}\subset U_{\phi_{l/k}(i),k}$, these groups are
equal if $i>H$, and $(U_{\phi_{l/k}(i),k}:N_{L/K}U_{i,l})=[l:k]=[L:K]$ if $i\le H$ because
the last residue fields are quasi-finite.

4. In the remaining case $L/K$ is infernal. Since both $K$ and $L$
are standard, $L=k\<t\>$, $K=k\<N_{L/K}t\>$. $L/K$ is ferociously
ramified, and Proposition 5.1.3 is applicable. This proves this
case of Proposition for $\al\ge(\c,\infty)$. For
$\al<(\c,\infty)$, notice that $N_{L/K}(\{t\}\cdot
U_{i,k})=\{N_{L/K}t\}\cdot U_{i,k}$.\qed
\enddemo

\subhead 6.4. Main theorems\endsubhead

\proclaim{Theorem 1}
Let $L/K$ be a finite abelian extension, $\al\in\BI_2$. Then $N_{L/K}Fil_\al K_2\tp L$
is a subgroup in $Fil_{\Phi_{L/K}(\al)}K_2\tp K$ of index $|Gal(L/K)_{\al}|$.
Furthermore,
$$
Fil_{\Phi_{L/K}(\al)}K_2\tp K\cap N_{L/K}K_2\tp L=N_{L/K}Fil_\al K_2\tp L.
$$
\endproclaim

\demo{Proof}
Use induction on $[L:K]$. If $[L:K]$ is
prime, then the assertion of Theorem is just Proposition 6.3.
In the general case let $h$ be the minimal ramification jump in $L/K$,
$G'=\cup_{\al>h}G_\al$, $G''\supset G'$ such that $(G:G')$ is prime, $K'=L^{G'}$.

Let $\al> h$. By assumption of induction, $N_{L/K'}Fil_\al K_2\tp L$
is a subgroup in $Fil_{\Phi_{L/K'}(\al)}K_2\tp K'$  of index $|\Gal(L/K')_{\al}|=|\Gal(L/K)_{\al}|$.
Applying $N_{K'/K}$, we obtain that $N_{L/K}Fil_\al K_2\tp L$ is a subgroup
of index $\le|\Gal(L/K)_{\al}|$ in
$$
Fil_{\Phi_{K'/K}(\Phi_{L/K'}(\al))}K_2\tp K=
Fil_{\Phi_{L/K}(\al)}K_2\tp K.
$$
If this index is less than $|\Gal(L/K)_{\al}|$, we easily obtain that
$(Fil_{\Phi_{L/K}(h)}K_2\tp K:N_{L/K}Fil_h K_2\tp L)<|\Gal(L/K)_{h}|\cdot [K':K]=[L:K]$.
On the other hand, also by induction on $[L:K]$, the norm map induces an epimorphism
of $K_2\tp L/Fil_h K_2\tp L$ onto $K_2\tp K/Fil_{\Phi_{L/K}(h)}K_2\tp K$, whence
$(K_2\tp K:N_{L/K}K_2\tp L)<[L:K]$, a contradiction with class field theory.

In the case $\al\le h$ we have $(Fil_{\Phi_{L/K}(\al)}K_2\tp K:N_{L/K}Fil_\al K_2\tp L)\ge[L:K]$
by class field theory, and
$$
\split
&(Fil_{\Phi_{L/K}(\al)}K_2\tp K:N_{L/K}Fil_\al K_2\tp L)\\
&\qquad=
(Fil_{\Phi_{L/K}(\al)}K_2\tp K:N_{K'/K}Fil_{\Phi_{L/K'}(\al)} K_2\tp K')\times\\
&\qquad\quad\times (N_{K'/K}Fil_{\Phi_{L/K'}(\al)}K_2\tp K':N_{L/K}Fil_\al K_2\tp L) \\
&\qquad\le [K':K][L:K']=[L:K].\qed
\endsplit
$$
\enddemo

\proclaim{Theorem 2}
Let $L/K$ be a finite abelian extension, $\Theta\:K_2\tp K/N_{L/K}K_2\tp L\to \Gal(L/K)$
the reciprocity map. Then
$$
\Theta(Fil_\al K_2\tp K \mod N_{L/K}K_2\tp L)= \Gal(L/K)^\al
$$
for any  $\al\in\BI_2$.
\endproclaim

\demo{Proof}
Let $G=\Gal(L/K)$, $\be<\al$, $K'=L^{G^\al}$.
Then $\Gal(K'/K)^\al$ is trivial, and
 $Fil_\al K_2\tp K\subset N_{K'/K}Fil_{\Psi_{K'/K}(\al)} K_2\tp K'$ by Theorem 1.
It follows
$$
\Theta(Fil_\al K_2\tp K \mod N_{L/K}K_2\tp L)\subset  \Gal(L/K')=\Gal(L/K)^\al.
$$
It remains to compare the indices of two subgroups by means of second assertion
in Theorem 1.\qed
\enddemo

\bigskip
The Department of Mathematics and Mechanics,
St. Petersburg University,
Bibliotechnaya pl. 2,
St. Petersburg 198904, Russia.
\bigskip
E-mail: igor\@zhukov.pdmi.ras.ru, igor\_zh\@hotmail.com

\Refs

\widestnumber\key {HLF}

\ref\key{A1}
\by V. A. Abrashkin 
\paper Towards explicit description of ramification filtration in the 2-dimen\-si\-onal case 
\paperinfo Preprint of Nottingham University  00-01 (2000)
\finalinfo to appear in Proceedings of the conference ``Ramification theory
of arithmetic schemes'' (Luminy, 1999)
\endref 
{\tt\ http://www.maths.dur.ac.uk/pure/ps/va-drf.ps}

\ref\key{A2}
\bysame
\paper Ramification theory for higher dimenional local fields
\paperinfo Preprint of Durham University (2002)
\endref
{\tt\ http://www.maths.dur.ac.uk/pure/ps/va-mmud1.ps}

\ref\key{E}\by H. Epp
\paper Eliminating wild ramification
\jour  Invent. Math.  \vol 19 \yr 1973  \pages 235--249    \endref

\ref\key  F1
\by I. B. Fesenko
\paper Class field theory of multidimensional local fields
of characteristic zero with residue field of positive characteristic
\jour Algebra i Analiz
\vol 3\yr1991     \issue3
\pages 165--196
\transl
\jour St. Petersburg Math.~J.
\vol3 \yr1992
\pages 649--678
\endref

\ref\key  F2
\bysame
\paper Abelian local $p$-class field theory
\jour Math. Ann.
\vol301\yr1995
\pages561--586
\endref

\ref\key{F3}\bysame
\paper Sequential topologies and quotients of Milnor  $K$-groups
of higher local fields
\jour Algebra i analiz
\vol13\yr2001
\pages198--228
\transl\nofrills{English translation to appear in}
\jour 
\jour St. Petersburg Math.~J.
\vol13 \yr2002
\issue3
\endref



\ref\key  FV
\by I. B. Fesenko, S. V. Vostokov
\book Local fields and their extensions: a constructive approach
\publ AMS
\publaddr Providence
\yr 1993
\endref



\ref\key H
\by O.~Hyodo
\paper Wild ramification in the imperfect residue field case
\jour Adv. Stud. Pure Math. \vol12\yr1987
\pages287--314
\endref

\ref\key HLF
\eds I.~B. Fesenko, M.~Kurihara
\book Invitation to higher local fields
\bookinfo Geometry and Topology Monographs, vol. 3
\yr2000
\endref

{\tt\ http://www.maths.warwick.ac.uk/gt/gtmcontents3.html}

\ref\key Kn
\by B. Kahn
\paper L'\,anneau de Milnor d'\,un corps local \`a corps residuel parfait
\jour Ann. Inst. Fourier, Grenoble
\vol34\yr1984\issue4\pages19--65
\endref

\ref\key{KZ}\by M. V. Koroteev, I. B. Zhukov \paper Elimination of
wild ramification \jour Algebra i
Analiz\vol11\issue6\yr1999\pages153--177 \transl\jour St.
Petersburg Math. J.\vol11\yr2000
\endref

\ref\key{L}\by V. G. Lomadze
\paper On the ramification theory of two dimensional local fields
\jour  Mat. Sb.  \vol 109 \yr 1979 \issue 3 \pages 378--394
\transl\jour Math. USSR-Sb. \vol 37 \yr1980
\endref



\ref
\key  MZ1
\by A. I. Madunts, I. B. Zhukov
\paper Multidimensional complete fields:
topology and other basic constructions
\jour Trudy S.-Peterb. Mat. Obsch.
\vol3 \yr1995
\pages4--46
\transl
\jour Amer. Math. Soc. Transl. (Ser.~2)
\vol165\yr1995\pages1--34
\endref

\ref
\key MZ2
\bysame
\paper Additive and multiplicative expansions
 in multidimensional local fields
 \jour Zapiski Nauchn. Sem. S.-Peterburg. Otdel.
  Mat. Inst. Steklov. (POMI)
  \vol272 \yr2000\pages186--196
  \endref

\ref\key  P
\by A. N. Parshin
\paper Local class field theory
\jour Trudy Mat. Inst. Steklov
\vol 165
\yr 1984
\pages 143--170
\transl
\jour Proc. Steklov Inst. Math.
\yr 1985 \vol 165\issue3 \pages 157--185
\endref

\ref\key S
\by J.-P. Serre
\book Corps Locaux
\bookinfo 2nd ed.
\publ Hermann
\publaddr Paris
\yr 1968
\endref


\ref\key Z1
\by I. B. Zhukov
\paper Structure theorems for complete fields
\jour Trudy S.-Peterb. Mat. Obschestva
\vol3\yr1995
\pages215--234
\transl
\jour Amer. Math. Soc. Transl. (Ser. 2)
\vol165\yr1995\pages175--192
\endref

\ref\key Z2
\bysame
\paper Milnor and topological $K$-groups of multidimensional complete fields
\jour Algebra i Analiz
\vol9\yr1997\issue1\pages98--147
\transl\jour St. Petersburg Math. J.
\vol9\yr1998\issue1\pages69--105
\endref


\endRefs
 \bye